\documentclass[reqno, twoside,a4paper]{amsart}
\usepackage{style}
\renewcommand{\P}{\mathcal{P}}
\newcommand{\ddd}{\stackrel{\rm def}{=}}

\begin{document}

\title[Orthogonal polynomials]{Zero sets, entropy, and pointwise asymptotics\\ of orthogonal polynomials}
\author{Roman Bessonov}
\author{Sergey Denisov}

\address{
\begin{flushleft}
Roman Bessonov: bessonov@pdmi.ras.ru\\\vspace{0.1cm}
St.\,Petersburg State University\\  
Universitetskaya nab. 7-9, 199034 St.\,Petersburg, RUSSIA\\
\vspace{0.1cm}
St.\,Petersburg Department of Steklov Mathematical Institute\\ Russian Academy of Sciences\\
Fontanka 27, 191023 St.Petersburg,  RUSSIA\\
\end{flushleft}
}
\address{
\begin{flushleft}
Sergey Denisov: denissov@wisc.edu\\\vspace{0.1cm}
University of Wisconsin--Madison\\  Department of Mathematics\\
480 Lincoln Dr., Madison, WI, 53706, USA\vspace{0.1cm}\\
\vspace{0.1cm}
Keldysh Institute of Applied Mathematics\\ Russian Academy of Sciences\\
Miusskaya pl. 4, 125047 Moscow, RUSSIA\\
\end{flushleft}
}

\thanks{
The work of RB done in Sections 2, 3 is supported by grant RScF 19-11-00058 of the Russian Science Foundation.
The work of SD done in Sections 1, 4 is supported by grant RScF-19-71-30004 of the Russian Science Foundation. His research conducted in the rest of the paper is supported by the grant NSF DMS-1764245 and by Van Vleck Professorship Research Award. }


\begin{abstract} Let $\mu$ be a measure from Szeg\H{o} class on the unit circle $\T$ and let $\{f_n\}$ be the family of Schur functions generated by $\mu$. In this paper, we prove a version of the classical Szeg\H{o}'s formula which
controls the oscillation of $f_n$ on $\T$ for all $n \ge 0$. Then, we focus on an analog of Lusin's conjecture for polynomials $\{\phi_n\}$ orthogonal with respect to  measure $\mu$ and prove that pointwise convergence of $\{|\phi_n|\}$ almost everywhere on $\T$ is equivalent to a certain condition on zeroes of $\phi_n$.  
\end{abstract}

\maketitle

\section{Introduction}
Consider a probability measure $\mu$ on the unit circle $\T = \{z \in \C: |z| = 1\}$ of the complex plane $\C$. The Schur function of $\mu$ is the analytic function $f$ in the open unit disk $\D = \{z \in \C: |z| < 1\}$ defined by the relation
\begin{equation}\label{eq00}
\frac{1 + zf(z)}{1-zf(z)} = \int_{\T}\frac{1 + \bar \xi z}{1 - \bar \xi z}\,d\mu(\xi), \qquad z \in \D.
\end{equation}
Taking the real part of both sides of \eqref{eq00} and using the Schwarz lemma, it is not difficult to see that $|f(z)| \le 1$ for all $z \in \D$. In particular, the function $f$ has non-tangential boundary values (to be denoted by the same letter $f$) almost everywhere on the unit circle $\T$. Set $f_0 = f$ and denote  the Schur iterates of $f$ by $f_n$: 
\begin{equation}\label{eq03}
zf_{n+1}(z) = \frac{f_n(z) - f_n(0)}{1 - \ov{f_n(0)}f_n(z)}, \qquad z \in \D, \qquad n \ge 0. 
\end{equation}
Schur's algorithm \eqref{eq03} produces an infinite family $\{f_n\}_{n \ge 0}$ of  analytic contractions unless $\mu$ is supported on a finite subset of $\T$, or, equivalently, $f$ is a finite Blaschke product. Knowing coefficients $f_k(0)$ for $0 \le k \le n$, one can set $f_{n+1} = 0$ and reverse the recursion in \eqref{eq03} to obtain an efficient approximation to $f$ in $\D$ by a rational contraction of degree $n$, see \cite{KH01}.

\medskip

Let $m$ be the Lebesgue measure on the unit circle $\T$ normalized by $m(\T) = 1$, and let $\mu = w\,dm + \mus$ be the decomposition of $\mu$ into the absolutely continuous and singular parts. The measure $\mu$ is said to belong to the Szeg\H{o} class $\szc$ if $\log w \in L^1(\T)$. To every measure $\mu \in \szc$, we associate the entropy function 
\begin{equation}\label{eq102}
\K(\mu, z) = \log\P(\mu, z) - \P(\log w, z), \qquad z \in \D,
\end{equation}
where $\P$ stands for the harmonic extension to $\D$:
$$
\P(\mu, z) = \int_{\T}\frac{1 - |z|^2}{|1 - \bar \xi z|^2}\,d\mu(\xi), 
$$
and we set $\P(v, z) = \P(v\,dm,z)$ for $v \in L^1(\T)$. Roughly speaking, $\K(\mu, z)$ measures  a ``size of oscillation'' of  $\mu$ on the arc $\{\xi \in \T: |\xi - a_z| \le 1-|z|\}$, $a_z = z/|z|$. By Jensen's inequality, we have $\K(\mu, z) \ge 0$ for every $z \in \D$ and $\K(\mu, z) = 0$ if and only if $\mu = m$. Notice also that $\K(\mu, \cdot )$ is superharmonic in~$\D$ and its nontangential boundary value is zero almost everywhere on $\T$.  

\medskip 

The celebrated Szeg\H{o} theorem says that a probability measure $\mu$ on the unit circle $\T$ belongs to the Szeg\H{o} class $\szc$ if and only if 
$\sum_{n \ge 0}|f_n(z)|^2<\infty$ 
for some (and then for every) $z \in \D$. Moreover, in the latter case we have
\begin{equation}\label{eq04}
\K(\mu, 0) = -\int_{\T}\log w\,dm =  -\log \prod_{n \ge 0}(1 - |f_n(0)|^2).
\end{equation}
This result has many equivalent reformulations, see, e.g., Section 2.7.8 in \cite{simonbook}. Our first aim is to extend formula \eqref{eq04} to the whole unit disk $\D$.
\begin{Thm}\label{t1}
Let $\mu \in \szc$ and let $\{f_n\}$ be the Schur family of $\mu$. Then
\begin{equation}\label{eq02}
\K(\mu, z) = \log \prod_{n \ge 0} \frac{1 - |zf_n(z)|^2}{1 - |f_n(z)|^2}, \qquad z \in \D. 
\end{equation}
\end{Thm}
Substituting $z = 0$ into \eqref{eq02}, we get \eqref{eq04}. As an immediate consequence of \eqref{eq02}, we see that $\sup_{n \ge 0}|f_n(z)|$ cannot be close to $1$ if $\K(\mu, z)$ is small. 

\medskip

Given a measure $\mu \in \szc$ and its Schur family $\{f_n\}$, we let $\mu_n$ denote the probability measure on $\T$ whose Schur function $f$ in \eqref{eq00} equals $f_n$. A standard problem in the field is to relate properties of $\mu_n$ to those of $\mu$ when $n$ is large. The following inequality is another immediate consequence of Theorem~\ref{t1}.
\begin{Cor}\label{c0}
We have $\K(\mu_n, z) \le \K(\mu,z)$ for all $n \ge 0$ and all $z \in \D$.
\end{Cor}
\noindent Indeed, due to Theorem~\ref{t1} and Schur's algorithm, we have 
$$
\K(\mu_n, z) = \log \prod_{k \ge n} \frac{1 - |zf_k(z)|^2}{1 - |f_k(z)|^2}.
$$
Since the terms in the product above are greater than $1$, we have $\K(\mu_n, z) \le \K(\mu, z)$.

\medskip

Theorem \ref{t1} implies a uniform bound for oscillation of Schur family generated by a Szeg\H{o} measure. 
\begin{Thm}\label{t3}
Suppose $\mu \in \szc$ and let $\{f_n\}$ be the family of Schur functions of~$\mu$. Then, we have
$$
\P(|f_n - f_n(z)|, z) \le c\sqrt{\K(\mu,z)}, \qquad z \in \D,
$$
with an absolute constant $c$ and all $n \ge 0$.
\end{Thm}
Let us now turn to an application of these results in the study of asymptotic behavior of orthogonal polynomials. To every measure $\mu \in \szc$ we associate the Szeg\H{o} function $D_\mu$. This is the outer function in $\D$ with modulus $\sqrt{w}$ on $\T$:
$$
D_{\mu}(z) = \exp\left(\int_{\T}\log\sqrt{w(\xi)} \cdot\frac{1+\bar\xi z}{1-\bar\xi z}\,dm(\xi)\right), \qquad z \in \D.
$$  
The family $\{\phi_n\}_{n \ge 0}$ of orthonormal polynomials in $L^2(\mu)$ is defined by 
\begin{equation}\label{sad3}
\deg \phi_n=n,\qquad k_n={\rm coeff}_n\phi_n>0, \qquad  (\phi_n, \phi_k)_{L^2(\mu)} = \delta_{n,k},
\end{equation}
where $\delta_{n,k}$ is the Kronecker symbol and ${\rm coeff}_jQ$ denotes the coefficient at the power $z^j$ in polynomial $Q$. Let also $\phi_n^*(z) = z^n \ov{\phi_n(1/\bar z)}$ denote the reversed orthogonal polynomial.  Due to a version of Szeg\H{o} theorem, we have $\mu \in \szc$ if and only if for some (and then for every) $z \in \D$ we have 
\begin{equation}\label{eq005}
\lim_{n \to +\infty} \phi_{n}^*(z) = D_{\mu}^{-1}(z).
\end{equation} 
A well-known conjecture in the theory of orthogonal polynomials on the unit circle (an analog of Lusin's conjecture for trigonometric series)  asks whether  \eqref{eq005} holds
for almost every $z \in \T$. As usual, for $z \in \T$ we understand $D_{\mu}^{-1}(z)$ as non-tangential boundary value. While not stated explicitly, the conjecture goes back to works of Bernstein, Szeg\H{o}, and Steklov who studied asymptotics of orthogonal polynomials. Recently, it attracted more attention due to its connection to  ``nonlinear Carleson problem'' in the scattering theory, see, e.g., \cite{ChK98}, \cite{ChK01}, \cite{MTT03}, \cite{OSTTW}. In the theorem below, we relate pointwise asymptotics of $\{\phi_n(z)\}$, $z\in \T$,  to the distribution of their zeroes near the unit circle. Our analysis is based on controlling oscillation of  Schur functions $\{f_n\}$ in terms of the entropy function $\K$ in \eqref{eq102}. The introduction of $\K$ was inspired by recent analysis of Szeg\H{o} condition for canonical systems \cite{BD2019}, \cite{BD2017}.

\medskip

Given a parameter $\rho\in (0,1)$ and a point $\xi\in \T$, define the Stolz angle $S_\rho^*(\xi)$ to be the convex hull of $\rho\D$ and $\xi$. Here is our main result.
\begin{Thm}\label{t2}
Let $\mu \in \szc$ and  $Z(\phi_n) = \{z \in \D: \; \phi_n(z) = 0\}$. Take any $a>0$ and denote $r_{a,n}= 1-a/n$. Then, for almost every $\xi\in \T$, the following assertions are equivalent:
\begin{itemize}
\item[$(a)$] $\lim_{n \to \infty} |\phi_{n}^{*}(\xi)|^2= |D_{\mu}^{-1}(\xi)|^2$,
\item[$(b)$] $\lim_{n \to \infty} \dist(Z(\phi_n), \xi)\,n = +\infty$,
\item[$(c)$] $\lim_{n \to \infty} f_n(r_{a,n}\xi) = 0$,
\item[$(d)$] $\lim_{n\to \infty}\sup_{z\in S_\rho^*(\xi)}|f_n(z)|=0$ for every $\rho\in (0,1)$.
\end{itemize} 
\end{Thm}
The paper is organized as follows. In Section \ref{s2}, we prove Theorem \ref{t1} and discuss its corollaries. Theorem \ref{t3} is proved in Section \ref{s3}.  In Section \ref{s4}, we collect some facts about finite sums of Poisson kernels that will be used in  Section \ref{s5} to prove Theorem~\ref{t2}.

\medskip

\section{Proof of Theorem \ref{t1} and some corollaries}\label{s2}
We start by giving an expression for $\K(\mu,z)$ in terms of $f$, the Schur function of measure $\mu$. 
\begin{Lem}\label{l0} If $\mu\in \szc$ and $f$ is its  Schur function, then
\begin{equation}\label{sad2}
\K(\mu, z) = \int_{\T}\log\left(\frac{1 - |zf(z)|^2}{1 - |f(\xi)|^2}\right)\,\frac{1-|z|^2}{|1 - \bar \xi z|^2}dm(\xi),
\end{equation}
for every $z \in \D$.
\end{Lem}
\beginpf Let $w$ be the density of $\mu$ with respect to $m$. Taking the real part of both sides of \eqref{eq00}, we obtain
$$
\frac{1-|zf(z)|^2}{|1 - z f(z)|^2} = \P(\mu,z), \qquad z \in \D.
$$
Hence, $w = \frac{1-|f|^2}{|1 - \xi f|^2}$ almost everywhere on $\T$. Then, the mean value formula for harmonic function $\log|1 - z f|^2$ implies
\begin{align*}
\K(\mu, z)
&=\log\frac{1-|zf(z)|^2}{|1 - z f(z)|^2} - \int_{\R}\log \frac{1-|f(\xi)|^2}{|1 - \xi f(\xi)|^2}\frac{1-|z|^2}{|1 - \bar\xi z|^2}\,dm(\xi)\\
&=\log (1-|zf(z)|^2)-\int_{\T}\log(1-|f(\xi)|^2)\,\frac{1-|z|^2}{|1 - \bar \xi z|^2}dm(\xi)\\
&=\int_{\T}\log\left(\frac{1 - |zf(z)|^2}{1 - |f(\xi)|^2}\right)\,\frac{1-|z|^2}{|1 - \bar \xi z|^2}dm(\xi).
\end{align*}
The lemma follows. \qed

\medskip

Now, let $\mu \in \szc$ and sequence $\{f_n\}_{n \ge 0}$ be the family of Schur functions generated by $\mu$ via the Schur's algorithm \eqref{eq03}. Denote by $\mu_n$ the probability measure on $\T$ whose Schur function coincides with $f_n$. Its existence follows if we notice that the function defined for $z\in \D$ by  
\begin{equation}\notag
\frac{1 - |zf_n(z)|^2}{|1 - zf_n(z)|^2} = \Re\left(\frac{1+zf_n(z)}{1-zf_n(z)}\right)
\end{equation}
is a nonnegative harmonic function in  $\D$ and therefore it is a Poisson integral of a unique nonnegative measure on $\T$. This is our $\mu_n$. Taking $z = 0$ in the  formula $\P(\mu_n,z) = \frac{1 - |zf_n(z)|^2}{|1 - zf_n(z)|^2}$, we get $\mu_n(\T) = 1$ so $\mu_n$ is a probability measure.

It is clear from construction that the Schur family of $\mu_n$ is $\{f_{n+k}\}_{k \ge 0}$.
After making these observations, we proceed with the proof of Theorem \ref{t1}.

\medskip

\noindent{\bf Proof of Theorem \ref{t1}.} For a measure $\mu \in \szc$, consider the family of Schur functions $\{f_n\}_{n \ge 0}$ and associated probability measures $\{\mu_n\}_{n \ge 0}$. By Szeg\H{o} theorem, we have $\sum_{k \ge 0}|f_{k}(0)|^2 < \infty$. It follows (again from the Szeg\H{o} theorem) that $\mu_n \in \szc$ and 
\begin{equation}\label{eq08}
\int_{\T}\log(1 - |f_n(\xi)|^2)\,dm = \log \prod_{k \ge n}(1 - |f_{k}(0)|^2) \to 0, \qquad n \to +\infty. 
\end{equation}  
In particular, functions $f_n$ tend to zero in Lebesgue measure on $\T$ and, since they are uniformly bounded, we have $\lim_{n \to \infty}f_{n}(z) = 0$ for every $z \in \D$. From \eqref{eq08} and Lemma \ref{l0}, we get  
$$
\K(\mu_n, z) = \int_{\T}\log \left(\frac{1 - |z f_n(z)|^2}{1-|f_n(\xi)|^2}\right)\frac{1-|z|^2}{|1 - \bar \xi z|^2}\,dm(\xi) \to 0, \qquad n \to +\infty,
$$
for every $z \in \D$. Thus, to prove Theorem \ref{t1}, we only need to check that 
\begin{equation}\label{eq0}
\K(\mu, z) = \K(\mu_{1}, z) + \log\frac{1-|zf(z)|^2}{1-|f(z)|^2}
\end{equation} 
and then iterate this formula.  From \eqref{sad2}, we have
\begin{align*}
\K(\mu, z) &= \log(1 - |zf(z)|^2) - \P(\log(1 - |f(\xi)|^2), z),\\
\K(\mu_1, z) &= \log(1 - |zf_1(z)|^2) - \P(\log(1 - |f_1(\xi)|^2), z),
\end{align*}
for every $z \in \D$.  Due to Schur's algorithm \eqref{eq03}, one can write
$$
zf_1(z) = \frac{f(z) - f(0)}{1 - \ov{f(0)}f(z)},
\qquad 
1-|zf_1(z)|^2 = \frac{(1-|f(0)|^2)(1 - |f(z)|^2)}{|1 - \ov{f(0)} f(z)|^2}.
$$
Using this computation, the mean value formula, and identity $|\xi| = 1, \xi\in\T$, we get
\begin{align*}
\K(\mu_1, z) 
&= \log\frac{(1-|f(0)|^2)(1 - |f(z)|^2)}{|1 - \ov{f(0)}f(z)|^2} - \P\left(\log\frac{(1-|f(0)|^2)(1 - |f(\xi)|^2)}{|1 - \ov{f(0)} f(\xi)|^2}, z\right)\\
&= \log(1 - |f(z)|^2) - \P\left(\log(1 - |f(\xi)|^2), z\right)\\
&= \log\frac{1 - |f(z)|^2}{1 - |zf(z)|^2} + \K(\mu,z),
\end{align*}
as required. \qed

\medskip

\begin{Cor}\label{c1}
Let $\mu \in \szc$ and $\{f_{n}\}_{n \ge 0}$ be the Schur family of $\mu$. Then, 
$$\K(1 - |f_n(\xi)|^2, z) \le \K(\mu, z)$$ 
for every $z \in \D$ and $n \ge 0$.
\end{Cor}
\beginpf 
Since $|zf_n(z)|^2$ is subharmonic in $\D$, we get
\[
\P(|\xi f_n(\xi)|^2,z)\ge |zf_n(z)|^2\,.
\]
Therefore,
\[
\log(1-|zf_n(z)|^2)\ge \log \P(1-|f_n(\xi)|^2,z)\,.
\]
So, applying Lemma \ref{l0} to measure $\mu_n$, we have
\begin{align*}
\K(\mu_n, z) 
&= \log(1 - |zf_n(z)|^2) - \int_{\T}\log(1 - |f_n(\xi)|^2)\,\frac{1-|z|^2}{|1 - \bar \xi z|^2}dm(\xi)\\
&\ge \log\P(1 - |f_n(\xi)|^2, z) - \P(\log(1 - |f_n(\xi)|^2), z)\\
&=\K(1 - |f_n(\xi)|^2, z).
\end{align*}
It remains to use Corollary \ref{c0}.  \qed

\medskip

Let $\alpha \in \T$, and let $f$ be the Schur function of a measure $\mu \in \szc$. Then, the family of measures $\mu_{\alpha}$ defined by
$$
\P(\mu_\alpha, z) = \Re\left(\frac{1+\alpha z f(z)}{1 - \alpha zf(z)}\right), \qquad z \in \D,
$$
is called the Aleksandrov-Clark family of $\mu$. From \eqref{eq00}, we see that $\alpha f$ is the Schur function of $\mu_\alpha$. 
\begin{Cor}\label{c3}
Let $\mu \in \szc$ and let $\{f_{n}\}_{n \ge 0}$ be the Schur family of $\mu$. Then, for every $z \in \D$, the entropy $\K(\mu, z)$ depends only on absolute value of $f(z)$. In particular, we have $\K(\mu, z) = \K(\mu_\alpha, z)$ for every $\alpha \in \T$.
\end{Cor} 
\beginpf This follows from \eqref{sad2}. \qed 

\medskip

The case $\alpha=-1$ in Corollary \ref{c3} corresponds to the ``dual measure'' $\mu_{\rm dual}$, playing an important role in the theory of orthogonal polynomials on the unit circle. The measure $\mu_{\rm dual}$ is defined by 
\[
\int_{\T}\frac{1 + \bar \xi z}{1 - \bar \xi z}\,d\mu_{\rm dual}(\xi)=\left(
\int_{\T}\frac{1 + \bar \xi z}{1 - \bar \xi z}\,d\mu(\xi)
\right)^{-1}, \qquad z \in \D.
\]
From \eqref{eq00}, we infer that the Schur function of $\mu_{\rm dual}$ equals $-f$. In particular, the last corollary yields
\begin{equation}\label{sd_1}
\K(\mu, z) = \K(\mu_{\rm dual}, z)\,, \qquad z \in \D.
\end{equation}

\medskip

It is well-known (see, e.g, Section 5 in \cite{KH01}) that orthonormal polynomials $\phi_n$ defined in \eqref{sad3} satisfy recurrence relations  
\begin{equation}\label{eq05}
\sqrt{1-|a_{n}|^2}\cdot \phi_{n+1}^* = \phi_{n}^* - za_{n}\phi_{n}, \qquad \phi_0 = \phi_0^* = 1, \qquad n \ge 0,  
\end{equation}
for coefficients $a_{n} = f_{n}(0)$ in $\D$, $\{f_n\}$ being the Schur family of $\mu$. Conversely, each sequence $\{a_k\}_{k \ge 0} \subset \D$ gives rise to a unique probability measure $\mu$ on $\T$ with infinite number of points in $\supp\mu$ such that its orthonormal polynomials satisfy relations \eqref{eq05}. In the next result we determine $\widehat\mu_{n,z}$, a variant of Bernstein-Szeg\H{o} approximation to $\mu$ such that $\K(\mu, z) = \K(\widehat\mu_{n,z}, z) + \K(\mu_{n+1}, z)$.

\medskip

\begin{Cor}
Let $n \ge 0$ and $z^\ast \in \D$. Consider the measure
$\widehat\mu_{n,z^\ast} = w_{n,z^\ast}\,dm$, where 
\begin{equation}\label{sad5}
w_{n,z^\ast}(\xi) = \frac{1 - |f_{n}(z^\ast)|^2}{|\phi_{n}^{*}(\xi) - \xi \overline{{f}_{n}(z^\ast)}\phi_{n}(\xi)|^2}, \qquad \xi \in \T.
\end{equation}  
Then, $\widehat\mu_{n,z^\ast}$ is a probability measure whose Schur  functions $\{\widehat f_{k}\}$ satisfy 
\begin{equation}\label{eq09}
\widehat{f}_{k}(z^\ast) = 
\begin{cases}f_k(z^\ast), & 0 \le k \le n,\\ 0, & k > n,\end{cases} 
\end{equation} 
at the point $z^\ast$. Moreover, we have $\K(\mu, z^\ast) = \K(\widehat\mu_{n,z^\ast}, z^\ast) + \K(\mu_{n+1}, z^\ast)$.
\end{Cor}
\beginpf Consider the family of orthonormal polynomials $\{\widehat\phi_j\}$ whose recurrence coefficients are given by 
$\widehat a_k = f_{k}(0)$ for $0 \le k \le n-1$, $\widehat a_{n} = f_n(z^*)$, and $\widehat a_{k} = 0$ for $k>n$. It is well-known that 
the measure $\nu = |\widehat\phi_{n+1}^*|^{-2}\,dm$ is a probability measure on $\T$ and its Schur functions $\{f_{\nu, k}\}_{k \ge 0}$ satisfy
 $f_{\nu, k}(0) = \widehat a_k$ for all $k \ge 0$. To see this, combine formulas $(4.17)$ and $(5.11)$ in \cite{KH01}. 
It follows that for all $w \in \D$ we have
$f_{\nu, n+1}(w) = 0$.  Therefore, from the definition of Schur's algorithm \eqref{eq03}, we have
\[
0=\frac{f_{\nu,n}(w)-f_{\nu,n}(0)}{1-\overline{ f_{\nu,n}(0)}f_{\nu,n}(w)}
\]
for all $w\in \D$ and so
$$ \qquad f_{\nu, n}(w) = f_{\nu, n}(0) = \widehat a_n = f_{n}(z^\ast).$$ 
Then, since $\widehat a_k = f_{k}(0)$ for all $0 \le k \le n-1$, we have 
$$f_{\nu, k}(z^*) = f_{k}(z^*), \qquad 0 \le k \le n-1$$ 
by Schur's algorithm \eqref{eq03} since $\{f_{\nu,k}\}$ and $\{f_k\}$ satisfy the same recursion at point $z^*$ when $k=0,1,\ldots,n-1$. To finish the proof, it remains to check that
$|\widehat\phi_{n+1}^*(\xi)|^{-2} = w_{n,z^*}(\xi)$ for $\xi \in \T$. To this end, observe that polynomials 
$\widehat\phi_{n}^*$ and $\phi_{n}^*$ are identical since the recurrence coefficients defining them are the same. Then, from \eqref{eq05} we get
$$
\sqrt{1-|\widehat a_{n}|^2}\cdot \widehat\phi_{n+1}^* = \widehat\phi_{n}^* - \xi\overline{\widehat a_{n}}\widehat\phi_{n}, \qquad  \widehat a_{n} = f_{n}(z^*),
$$
and \eqref{sad5} follows. \qed

\medskip

According to a theorem by Khrushchev (Theorem 3 in \cite{KH01}), the Schur function of the probability measure $|\phi_n^*|^2\,d\mu$ is equal to $b_n f_n$, where $b_n = \phi_n/\phi_n^*$ is the Blaschke product of order $n$. In other words, we have (formula (2.14) in \cite{KH01})
\begin{equation}\label{eq11}
\int_{\T}\frac{1+\bar\xi z}{1-\bar\xi z}|\phi_n^*(\xi)|^2\,d\mu(\xi) = \frac{1+ z b_n(z) f_n(z)}{1- z b_n(z) f_n(z)}, \qquad z \in \D,
\end{equation}
and hence (formula (1.18) in \cite{KH01})
\begin{equation}\label{sd_6}
|\phi_n^*(\xi)|^{2}w(\xi) = \frac{1 - |f_n(\xi)|^2}{|1 - \xi b_n(\xi) f_n(\xi)|^2}, \qquad \xi \in \T.
\end{equation}
Identity $|b_n(\xi)| = 1,\xi\in\T$ implies the following corollary.
\begin{Cor}
We have 
$$\K(|\phi_n^*(\xi)|^2\,d\mu, z) = \K(\mu_n, z) + \log\left(\frac{1 - |z b_n(z) f_n(z)|^2}{1 - |z f_n(z)|^2}\right)$$ 
for every $n \ge 0$ and $z \in \D$. 
\end{Cor} 
\beginpf Fix $n \ge 0$ and $z \in \D$. It follows from \eqref{sad2}  that 
\begin{align*}
\K(|\phi_n^*|^2\,d\mu, z) 
&= \log(1-|z b_n(z) f_n(z)|^2) - \P(\log(1-|b_n f_n|^2), z)\\
&= \log(1-|z b_n(z) f_n(z)|^2) - \P(\log(1-|f_n|^2), z)\\
&= \K(\mu_n, z) + \log\left(\frac{1 - |z b_n(z) f_n(z)|^2}{1 - |z f_n(z)|^2}\right),
\end{align*}  
as required. \qed

\medskip

Let us now consider the case when $\mu$ is absolute continuous and its density does not oscillate too much. We say that $w\in A_\infty^P(\T)$  if
\begin{equation}\label{eq10}
[w]_{\infty,P}=\sup_{z\in \D}\P(w,z)\exp\Bigl(-\P(\log w,z)\Bigr)<\infty\,.
\end{equation}
It is known that $A_\infty^P(\T)\subsetneq A_\infty(\T)$, where $A_\infty(\T)$ is the usual Muckenhoupt class. 
\begin{Lem} We have $w\in A_\infty^P(\T)$ if and only if $\sup_{z\in \D}\K(w\,dm,z)<\infty$. Moreover, the dual measure of $wdm$ is absolutely continuous -- i.e., $(w\,dm)_{\rm dual}=w_{\rm dual}dm$, and its density satisfies $w_{\rm dual}\in A_\infty^P(\T)$.
\end{Lem}
\beginpf The first statement is immediate from the definition. To prove the second one, we use \eqref{sd_1} and notice that  $\mu \in \szc$ and $\K(\mu,z)\in L^\infty(\D)$ imply that $\mu$ has no singular part and $\mu=wdm$ with $w\in A_{\infty}^P(\T)$. Indeed, if $\mu=w\,dm+\mus$ where  $\mus$ is the singular measure, then  
\[
\log\left( \P(\mus,z)+\P(w,z)\right)-\P(\log w,z)\le C, \quad z\in \D,
\]
by our assumptions.
This implies 
\[
\P(\mus,z)\le \P(\mus,z)+\P(w,z)\le C\exp\left(\P(\log w,z)\right)\le C\P(w,z),
\]
by Jensen inequality, hence, $\mus=0$. \qed
\begin{Cor} Let the probability measure $\mu$ be defined by $\mu=wdm$ and $w$ satisfy $w \in A_\infty^P(\T)$. If $\{f_n\}$ denotes the Schur family of $\mu$, then $1-|f_n|^2 \in A_{\infty}^P(\T)$ and $[1-|f_n|^2]_{\infty,P} \le [w]_{\infty,P}$, $n \ge 0$. 
\end{Cor}
\beginpf By Corollary \ref{c1}, for each $n \ge 0$ and $z \in \D$, we have 
$$
\log\P(1-|f_n|^2,z) - \P(\log (1-|f_n|^2),z) \le \K(\mu,z) \le \log [w]_{\infty,P}.
$$ 
It follows that $\log [1-|f_n|^2]_{\infty,P} \le \log [w]_{\infty,P}$. \qed

\medskip

\section{The space $\bmo_\eta$ and proof of Theorem \ref{t3}}\label{s3}
Given a function $\eta: \D \to [0, +\infty]$, we define the space $\bmo_\eta$ to be the set of functions $v \in L^1(\T)$ such that the following characteristic
$$
\|v\|_{\eta}^{*} = \inf\{c \ge 0: \P(|v - \P(v,z)|, z) \le c\eta(z), \; z \in \D\}
$$
is finite. The next result is a direct analogue of an estimate by M.~Korey (see Section 3.2 in \cite{Korey}).
\begin{Lem}\label{l07} Suppose that $v, e^v \in L^1(\T)$ and let $\P(e^v, z)/e^{\P(v, z)} = 1+\gamma$ for some $\gamma\ge 0$ and $z \in \D$. Then,
$$
\P(|v - \P(v,z)|, z) \le c 
\begin{cases}
\sqrt\gamma, & \gamma<1,\\
\log(1+\gamma), &\gamma\ge 1,
\end{cases}
$$
for an absolute constant $c$.
\end{Lem}
\beginpf The proof is an adaptation of the original argument in \cite{Korey}. For the reader's convenience, we reproduce it here. It suffices to prove the inequality 
$$
\P(|v - m_z(v)|, z) \le c 
\begin{cases}
\sqrt\gamma, & \gamma<1,\\
\log(1+\gamma), &\gamma\ge 1,
\end{cases}
$$ 
where $m_z(v)$ is the median value of $v$ on $\T$ with respect to the probability measure $\nu = (1 - |z|^2)/|1 - \bar \xi z|^2\,dm$. Adding a constant to $v$ if needed, one can assume that $m_z(v) = 0$. Then, there are two disjoint measurable subsets $E \subseteq \{\xi: v(\xi) \ge 0\}$ and $F\subseteq  \{\xi: v(\xi) \le 0\}$ of $\T$ such that $\nu(E) = \nu(F) = 1/2$. Set 
\begin{align*}
a = 2\P(\chi_E e^v, z), \quad b = 2\P(\chi_F e^v, z), \qquad a' = e^{2\P(\chi_E v, z)}, \quad b' = e^{2\P(\chi_F v, z)}. 
\end{align*}
By construction and by Jensen's inequality, one gets
$$
1+\gamma = \frac{\P(e^v, z)}{e^{\P(v, z)}} = \frac{a+b}{2\sqrt{a'b'}} \ge \frac{a'+b'}{2\sqrt{a'b'}}\,,
$$
which implies $a'/b' \le 1+\widetilde c\max(\sqrt{\gamma},\gamma^2)$ with an absolute constant $\widetilde c$. On the other hand, we have
$a'/b' = e^{2\P(\chi_E v, z) - 2\P(\chi_F v, z)}$. It follows that 
$$
\P(|v|, z) = \P(\chi_E v, z) - \P(\chi_F v, z) \le c\log(1+\max(\sqrt{\gamma},\gamma)),
$$ 
for another absolute constant $c$, as claimed. \qed

\medskip

Given a measure $\mu \in \szc$, we introduce the function 
\begin{equation}\label{sd_a1}
\eta(z) =\max\left(\sqrt{\K(\mu,z)}, \K(\mu,z)e^{\K(\mu,z)/2}\right),
\end{equation}
on the unit disk $\D$. The next lemma is crucial for later analysis.  
\begin{Lem}\label{l7}
Consider $\mu \in \szc$. Let $\{f_n\}$ be the Schur family of $\mu$ and $\{\phi_n\}$ be  orthogonal polynomials generated by $\mu$. Then the functions $\log|\phi_n^* - \xi f_{n}\phi_n|^2$ and $f_n$ belong to $\bmo_{\eta}$ for all $n \ge 0$ and
$$
\Bigl\|\log|\phi_n^* - \xi f_{n}\phi_n|^2\Bigr\|_{\eta}^{*} \le c, \quad \|f_n\|_{\eta}^{*} \le c, \qquad n \ge 0,
$$
with an absolute constant $c$.
\end{Lem}
\beginpf Consider the weight $v_n = 1 - |f_n|^2$ on the unit circle $\T$. By Corollary \ref{c1}, we have 
$\K(v_n, z) \le \K(\mu, z)$. Hence, applying  Lemma \ref{l07} to $v = \log v_n$ one has
$$
\P(|\log v_n - \P(\log v_n,z)|, z) \le c \max(\sqrt{\K(\mu,z)},\K(\mu,z))\le c  \eta(z), \qquad z \in \D.
$$
It follows that $\|\log v_n\|_{\eta}^{*} \le c$ for all $n \ge 0$. In a similar way, we get $\log w \in \bmo_{\eta}$. 
We now use \eqref{sd_6} to write
$$
\log w = \log v_n - \log|\phi_n^* - \xi f_{n}\phi_n|^2,
$$ 
hence $\log|\phi_n^* - \xi f_{n}\phi_n|^2 \in \bmo_{\eta}$ with the characteristic $\|.\|_{\eta}^{*}$ at most $2c$. Next, we use Jensen's inequality to write
$$
\P(\log(1 - |f_n|^2), z) \le \log \P(1 - |f_n|^2, z) = \log(1 - \P(|f_n|^2, z)).
$$
Therefore, applying Lemma \ref{l0} to measure $\mu_n$, one has
$$
\K(\mu_n, z) 
= \log(1 - |zf_n(z)|^2) - \P(\log(1 - |f_n|^2), z)
\ge \log\frac{1 - |f_n(z)|^2}{1 - \P(|f_n|^2, z)}.
$$
Since $\K(\mu_n, z) \le \K(\mu, z)$ by Corollary \ref{c0}, we have
$$
1 - |f_n(z)|^2 \le e^{\K(\mu, z)}(1 - \P(|f_n|^2, z)),
$$
which can be rewritten as
\[
e^{\K(\mu,z)}\P(|f_n|^2, z)-|f_n(z)|^2\le e^{\K(\mu,z)}-1\,.
\]
Since $\K\ge 0$, the following inequality holds
\[
\P(|f_n|^2, z)-|f_n(z)|^2\le e^{\K(\mu,z)}\P(|f_n|^2, z)-|f_n(z)|^2\le e^{\K(\mu,z)}-1\,.
\]
The last bound along with mean value formula for harmonic functions imply
\begin{align*}
\P(|f_n - \P(f_n,z)|^2, z) 
&=\P(|f_n - f_n(z)|^2, z),\\
&=\P(|f_n|^2, z) +|f_n(z)|^2-2\P(\Re (f_n \overline{f_n(z)}),z),\\
&=\P(|f_n|^2,z)-|f_n(z)|^2\le e^{\K(\mu, z)} - 1.
\end{align*}
By Cauchy-Schwarz inequality, 
we get
\[
\P(|f_n - \P(f_n,z)|, z)\le \sqrt{e^{\K(\mu, z)} - 1}\le c\eta(z)\,.
\]
That finishes the proof.
 \qed

\medskip

\noindent{\bf Proof of Theorem \ref{t3}.} By Lemma \ref{l7}, for every $n \ge 0$ we have 
$$
\P(|f_n - f_n(z)|, z) \le c\eta(z), \qquad z \in \D.
$$
On the other hand,  $\P(|f_n - f_n(z)|, z) \le 2$ since $|f_n| \le 1$ on $\D \cup \T$. This yields the statement of the theorem. \qed
 
\medskip

Next, we will estimate the harmonic conjugates of functions in $\bmo_\eta$.  Some notation is needed first. We denote $|E| = m(E)$ for Borel subsets of $\T$. If $I \subset \T$ is an arc with center at $\xi$, set $z_I=\xi(1-|I|)$ and denote by $2I$ the arc with center at $\xi$ such that $|2I|=2|I|$. We also let $\langle f\rangle_{I,P}=\P(f,z_I)$. For $u \in L^1(\T)$, we define the harmonic conjugate function $v$ by the formula
$$
v(\xi) = (Qu)(\xi) = \lim_{r \to 1}\int_{\T}u(\zeta)Q_{r}(\zeta, \xi)\,dm(\zeta), \quad 
Q_{r}(\zeta, \xi) = \Im \frac{1 + r\bar \zeta \xi}{1 - r\bar \zeta \xi}, \quad \xi\in \T\,.
$$
From the standard estimates for singular integrals, one knows that the limit exists almost everywhere on $\T$ and defines the function $Qu\in L^{1,\infty}(\T)$. Notice that the harmonic conjugate of a constant function is identically zero.  Finally, given real-valued $u\in L^1(\T)$, the function $u+i(Qu)$ is the nontangential boundary value of the function 
\[
\mathcal{F}(z)=\int_{\T}u(\zeta) \frac{1 + \bar \zeta z}{1 - \bar \zeta z}   \,dm(\zeta)
\]
analytic in $\D$. Function $\Re \mathcal{F}$ is Poisson extension of $u$ and harmonic conjugate of $u$ is the boundary value of $\Im \mathcal{F}$. Next,  we recall that, given a parameter $\rho\in (0,1)$, the symbol $S^*_\rho(\xi)$ denotes the convex hull of $\rho\D$ and a point $\xi\in \T$.

\medskip

Below we write $A \lesssim B$ for quantities $A$, $B$ if there is an absolute constant $c$ such that $A \le c B$. Notation $A \sim B$ is used when $A \lesssim B$ and $B \lesssim A$.

\begin{Lem}\label{l12}
Let $u \in \bmo_\eta$ and let $v$ be the harmonic conjugate of $u$. Let $I$ be an arc with center at $\xi_0 \in \T$. Then, there is a constant $c_I$  such that
$$
\frac{|\{\xi \in I: |v(\xi) - c_I| > t\}|}{|I|} \lesssim t^{-1} \|u\|^{*}_{\eta} \sum\nolimits_{j\ge 0}2^{-j}\eta(z_j),
$$
for some $z_j \in S_{0.9}^{*}(\xi_0)$ such that $|z_j - \xi_0| \sim 2^{j} |I|$, $j \ge 0$.
\end{Lem}

\beginpf Write $u = u_1 + u_2 + \langle u\rangle_{2I,P}$ for
$u_1 = \chi_{2I}(u - \langle u\rangle_{2I,P})$, $u_2 = \chi_{\T\setminus 2I}(u - \langle u\rangle_{2I,P})$, and denote by $v_1$, $v_2$ the harmonic conjugates of $u_1$, $u_2$, respectively.  
Since $Q$ is the continuous operator from $L^1(\T)$ to  $L^{1, \infty}(\T)$, we have 
\begin{align*}
m(\{\xi \in I: |v_1(\xi)| > t\}) 
&\lesssim t^{-1}\|u_1\|_{L^1(\T)}, \\
&\lesssim t^{-1}\langle|u - \langle u\rangle_{2I,P}|\rangle_{2I,P}|I|, \\
&\lesssim  t^{-1} \|u\|^{*}_{\eta}\cdot \eta(z_0) \cdot |I|,
\end{align*}
for  $z_0 =z_{2I}= (1 - 2|I|)\xi_0$. Next, we estimate the distribution function of~$v_2$. Put $c_I = \int_{\T\setminus 2I}u_2(\zeta)Q(\zeta, \xi_0)\,dm(\zeta)$
and  write for $\xi \in I$:
$$
v_2(\xi) - c_I = (Q u_2)(\xi) - c_I = \int_{\T\setminus 2I}u_2(\zeta)\Bigl(Q(\zeta, \xi) - Q(\zeta, \xi_0)\Bigr)\,dm(\zeta).
$$
Let us estimate the norm of $v_2 - c_I$ in $L^1(I)$ to later use Chebyshev inequality. For $k \ge 1$, denote by $I_k$ the arcs of $\T$ of size $2^{k}|I|$ with center at $\xi_0$. Then, for $\xi \in I$, $\zeta \in I_{k+1} \setminus I_k$, we have 
$$
|Q(\zeta, \xi) - Q(\zeta, \xi_0)| \lesssim \frac{|\xi - \xi_0|}{|(\zeta - \xi)(\zeta - \xi_0)|} \lesssim   \frac{|I|}{|I_k|^2}= \frac{1}{2^{2k}|I|}.
$$
Using this relation, we get
\begin{align*}
\int_{I}|v_2 - c_I|\,dm 
&\lesssim \sum_{k \ge 1} \frac{1}{2^{2k}|I|}\int_{I}\int_{I_{k+1}\setminus I_k}|u_2(\zeta)|\,dm(\zeta)\,dm(\xi)\\
&\lesssim \sum_{k \ge 1} \frac{1}{2^{2k}}\int_{I_{k+1}\setminus I_k}|u_2|\,dm.
\end{align*} 
Set $J_0 = 2I$ and let $J_k$, $k \ge 1$, be one of two arcs of $I_{k+1}\setminus I_k$ such that 
$$
\int_{I_{k+1}\setminus I_k}|u_2|\,dm \le 2 \int_{J_k}|u_2|\,dm.
$$
We have
\begin{align*}
\int_{J_k}|u_2|\,dm 
&\le |J_k| \cdot \Bigl(\langle |u - \langle u\rangle_{J_k,P}|\rangle_{J_k,P} + |\langle u\rangle_{J_k,P} - \langle u\rangle_{2I,P}|\Bigr) \\
&\lesssim  2^{k}|I| \cdot \Bigl(\|u\|_{\eta}^{*}\eta(z_k) + \sum_{j=1}^{k}|\langle u\rangle_{J_j,P} - \langle u\rangle_{J_{j-1},P}|\Bigr), 
\end{align*} 
where $z_k = (1 - |J_k|)\xi_k$ and $\xi_k$ denotes the center of $J_k$. Since $|\zeta - z_j|\sim |\zeta - z_{j+1}|$ for $\zeta\in \T$, we can write
$$
|\langle u\rangle_{J_j,P} - \langle u\rangle_{J_{j-1},P}|=|\P(u-\langle u\rangle_{J_{j-1},P}, z_j)| \lesssim
\P(|u-\langle u\rangle_{J_{j-1},P}|, z_{j-1}) \lesssim \|u\|_{\eta}^{*}\eta(z_{j-1})\,.
$$
Hence,
$$
\int_{J_k}|u_2|\,dm \lesssim  2^{k}\cdot|I|\cdot\|u\|_{\eta}^{*}\sum_{j=0}^{k}\eta(z_j).
$$
It follows that
$$
\frac{1}{|I|}\int_{I}|v_2 - c_I|\,dm \lesssim \|u\|_{\eta}^{*}\cdot\sum_{k \ge 1}2^{-k}\sum_{j=0}^{k}\eta(z_j) \lesssim 
\|u\|_{\eta}^{*}\cdot\sum_{j \ge 0}2^{-j}\eta(z_j).
$$
Now we collect estimates to get the bound
\begin{align*}
|\{\xi \in I: |v(\xi) - c_I| > 2t\}| 
&\le|\{\xi \in I: |v_1(\xi)| > t\}| + |\{\xi \in I: |v_2(\xi) - c_I| > t\}|\\
&\lesssim   t^{-1} \|u\|^{*}_{\eta}\cdot \eta(z_0) \cdot |I| + t^{-1}\int_{I}|v_2 - c_I|\,dm\\
&\lesssim  t^{-1} \|u\|^{*}_{\eta}\cdot |I| \sum\nolimits_{j\ge 0}2^{-j}\eta(z_j).
\end{align*}
The simple geometric considerations yield $z_j\in S_{0.9}^*(\xi_0)$ and the lemma is proved. \qed 

\medskip

\section{Sums of Poisson kernels}\label{s4}
In this section, we study the properties of finite sums of Poisson kernels. They will be used in the proof of Theorem \ref{t2}.  

\medskip

We denote by $C[a,b]$ the space of functions continuous on $[a,b]$. The following elementary result is well-known.
\begin{Lem}\label{sd_111}
Suppose the sequence $\{g_n\}$ of non-decreasing functions converges to a function $g\in C[a,b]$  on a dense subset of $[a,b]$. Then, $\{g_n\}$ converges to $g$ uniformly on $[a,b]$.
\end{Lem}

We start with the calculation which reveals the connection between the zeroes of the polynomial $\phi_n$ and the sum of Poisson kernels. Consider $b_n=\phi_n/\phi_n^*$. It can be written as
\[
b_n(z)=\alpha_nz^{l_n}\prod_{j=1}^{m_n} \frac{z-z_{j,n}}{z-\bar z_{j,n}^{-1}}, \quad \alpha_n>0,\quad l_n+m_n=n\,,
\]
where $\{z_{j,n}\}$ are zeroes of $\phi_n$ different from $0$. This is the product of M\"obius transforms each of which has an argument which is increasing monotonically on $\T$ since this transform is a conformal map of $\D$ onto $\D$. Calculating the derivative of its argument 
\begin{equation}\label{sd_a21}
\partial_\theta \arg b_n(e^{i\theta})=l_n+\Im \partial_\theta\left(\sum_{j=1}^{m_n} \log\left(\frac{e^{i\theta}-z_{j,n}}{e^{i\theta}-\bar z_{j,n}^{-1}}\right)\right)=l_n+\sum_{j=1}^{m_n} \frac{1-|z_{j,n}|^2}{|e^{i\theta}-z_{j,n}|^2}\,
\end{equation}
one can recognize the Poisson kernel as terms in the last sum. 
\begin{Lem}\label{l8}
Assume that $h_n$ are smooth functions on $(-\pi n,\pi n)$ with derivatives $h_n'$ given by
\begin{equation}\label{sd_7}
h'_n(t) = \frac{1}{n}\sum_{k=1}^{n}\frac{1-|z_{k,n}|^2}{|e^{it/n}-z_{k,n}|^2}, \qquad z_{k,n} \in \D. 
\end{equation}
If $\{h_n\}$ converges to a smooth function $h$ uniformly on compact subsets of $\R$, then $\{h'_n\}$ converges to $h'$ uniformly on compact subsets of $\R$.
\end{Lem}
\beginpf We will assume that the points $z_{k,n}$ are enumerated so that 
$$
|1-z_{k,n}| \le |1 - z_{k+1,n}|, \qquad 1 \le k \le n.
$$ 
Take an arbitrary $b>0$ and let $F_n=h_n'$. It suffices to show that $\{F_n\}$ converges to $h'$ uniformly over $[-b/2,b/2]$.
We write $F_n$ as $F_n=G_n+H_n$, where $G_n$ is the sum which corresponds to all terms (if any) for which $n|1-z_{k,n}|>1.9b$ and, respectively, terms in $H_n$ satisfy  $n|1-z_{k,n}|\le 1.9b$. For $G_n$, we have 
\begin{equation}
|G_n'(t)|\lesssim b^{-1}G_n(t)\label{sd_pa}
\end{equation} when $t\in [-b,b]$. Indeed, if $t\in [-b,b]$ and $n|1-z|\ge 1.9b$, we have
\begin{equation}\label{sd_22}
\partial_t\left(\frac{1}{|e^{it/n}-z|^2}\right)
=\frac{2}{n}\frac{\Re (i\bar{z}e^{it/n})}{|e^{it/n}-z|^4}
\le \frac{1}{n|e^{it/n}-z|^3} \lesssim \frac{1}{b|e^{it/n}-z|^2},
\end{equation}
which yields the required estimate. It follows that
$$
\limsup_n\int_{-b}^{b}|G_n'|\,dt \lesssim \limsup_n\int_{-2b}^{2b}b^{-1}G_n\,dt \le b^{-1}(h(2b)-h(-2b)).
$$ 
Thus, functions $\{G_n\}$ are uniformly bounded. The estimate \eqref{sd_pa} then implies that the set $\{G_n\}$ is also   equicontinuous on $[-b,b]$. Choose a subsequence $\{G_{n_j}\}$ which converges to some continuous function $G$ uniformly over $[-b,b]$. Then, $\{\int_{-b}^{x} G_{n_j}\,dt\}$ converges to $\int_{-b}^x G\,dt$ uniformly over $[-b,b]$ as well. Since we know by conditions of the lemma that
$\{\int_{-b}^x F_n\,dt\}$ converges uniformly to a smooth function, the sequence $\{\int_{-b}^x H_{n_j}\,dt\}$ also converges uniformly to a function continuous on $[-b,b]$.

\smallskip

Now let $z_{k,n_j}$, $k=1,\ldots,c(n_j)$ be all points that satisfy $n_j|1-z_{k,n_j}|\le 1.9b$. 
For every $z \in \D$ such that $n_j|1-z|\le 1.9b$, we have 
\[
\frac{1}{n_j}\int_{-2b}^{2b}\frac{1-|z|^2}{|e^{it/n_j}-z|^2}dt=\int_{-2b/n_j}^{2b/n_j}\frac{1-|z|^2}{|e^{i\tau}-z|^2}d\tau \ge c_b > 0,
\]
where the constant $c_b$ depends only on $b$.  It follows that
$$
\limsup_j c(n_j) \le c_{b}^{-1} \limsup_{j}\int_{-2b}^{2b} H_{n_j}\, dt\le (h(2b)-h(-2b))/c_b.
$$ 
Hence, $\limsup_j c(n_j) = N_b$ for some $N_b \ge 0$ (we set $N_b = 0$ if there are no zeroes $z_{k, n_j}$ such that $n_j|1-z_{k,n_j}|\le 1.9b$ for all $j$ large enough). Choosing a subsequence, on can assume that $c(n_j) = N_b$ for all $j$. If $N_b>0$, we set $\xi_{k,n_j}=in_j(1-z_{k,n_j})$ for every $k=1,\ldots, N_b$. Note that $\xi_{k,n_j}$ belong to $\C^+ = \{z \in \C: \; \Im z >0\}$ and, moreover, $|\xi_{k,n_j}|<1.9b$. Choosing again a subsequence, we may assume that $\{\xi_{k,{n_j}}\}$ converges to $\xi_k\in \overline{\C^+}$, $k=1,\ldots, N_b$. We claim that none of these limiting points belongs to the segment $[-b/2,b/2]$ on the real line. Indeed, if  $\xi_k\in [-b/2,b/2]$, then the sequence of functions
\[
\frac{1}{n_j}\frac{1-|z_{k,n_j}|^2}{|e^{it/n_j}-z_{k,n_j}|^2}, \qquad j=1,2,\ldots
\]
converges to $2\pi \delta_{\xi_k}$ in the weak-$\ast$ sense, where $c$ is a nonzero constant, see \eqref{sd_8}. This contradicts the fact that $\int_{-b}^{x} H_{n_j}dt$ converges uniformly to a continuous function on $[-b,b]$. Knowing that all limiting points $\xi_k$ are separated from the real line, it is easy to see that $\{H_{n_j}\}$ converges uniformly over $[-b/2,b/2]$. Thus, we can guarantee that some subsequence $\{F_{n_j}\}$ converges uniformly on $[-b/2, b/2]$. Denote its limit by $F$. Since 
\[
\lim_{n\to\infty}\int_{-b/2}^x F_n(t)\,dt = h(x)-h(-b/2)
\]
by assumption of the lemma, we get $F=h'$. In the standard way, we now get $\lim_{n \to +\infty} F_n = h'$ uniformly on compacts in $\R$ over the whole original sequence. Indeed, if this is not true, then there is $\eps>0$, $b>0$, $t_n\in [-b/2,b/2]$ and a sequence $\{m_n\}$ such that $|F_{m_n}(t_n)-h'(t_n)|>\eps$. However, by the argument above we can take a subsequence of $\{m_n\}$, call it $\{m_n'\}$, so that $\{F_{m_n'}\}$ converges to $h'$ uniformly on $[-b/2,b/2]$. This gives a contradiction. The lemma is proved. \qed

\medskip

\begin{Lem} \label{sd_12} 
In the previous lemma, suppose that $h(t)=t+c$ for all $t \in \R$ and some constant $c$. Then, 
$$
\lim\limits_{n\to\infty}\left(n\min \limits_{1 \le k \le n}|1-z_{k,n}|\right)=\infty.
$$
\end{Lem}
\beginpf 
Given $b\in \mathbb{N}$, let $\{z_{j,n}\}$, $j=1,\ldots,c(n,b)$, be all zeroes of $\phi_n$, counting multiplicity, that satisfy $n|1-z_{j,n}|<1.9b$, and set $c(n,b) = 0$ if there are no such zeroes. From the previous proof, we know that $\limsup_n c(n,b)<\infty$. We need to show that $\limsup_{n\to\infty} c(n,b)=0$ for every $b$. Suppose this is not the case and there is some $\widehat b$ such that $\limsup_{n\to\infty} c(n,\widehat b)\ge 1$. Then, there is a subsequence $\{n_k^{(1)}\}$ such that each $\phi_{n_k^{(1)}}(z)$ has exactly $c_1$ zeroes, counting multiplicity, at points $\{z_{j,n_k^{(1)}}\}, j=1,\ldots,c_1$, $c_1 \ge 1$, and all of these zeroes are inside the open disc of radius $1.9\widehat b/n$ centered at $1$. Using compactness argument, we can find $\{\widehat n_k^{(1)}\}$, a subsequence of  $\{n_k^{(1)}\}$, such that 
\[
\lim_{k\to\infty}\xi_{j,\widehat n_k^{(1)}} =\xi_j, \quad  \xi_{j,n}\ddd in(1-z_{j,n}),  \qquad j=1,\ldots,c_1\,
\]
and $c_1$ points $\{\xi_j\}$, counting multiplicity, all belong to the set $\{\xi\in \C^+: |\xi|\le 1.9\widehat b\}$. Notice that none of these points can be on the real line, this follows from the proof of the previous lemma. 

\medskip

Next, we look at zeroes of a polynomial $\phi_{\widehat n_k^{(1)}}(z)$ that belong to the annulus 
$$
1.9\widehat b\le n|1-z|<1.9(\widehat b+1).
$$ 
Applying the same argument, we can find $\{\widehat n_k^{(2)}\}$, a subsequence of $\{\widehat n_k^{(1)}\}$, for which each $\phi_{\widehat n_k^{(2)}}(z)$\
has exactly $c_2$ zeroes $\{z_{j,n_k^{(2)}}\}, j\in \{c_1,\ldots,c_1+c_2\}$ in that annulus and they all satisfy
 
 \[
\lim_{k\to\infty}\xi_{j,\widehat n_k^{(2)}} =\xi_j,   \qquad j=c_1,\ldots,c_1+c_2\,.
\]
Notice that it might be that $c_2=0$ but we always have $c_2<\infty$.

\medskip

We continue this process and find a subsequence $\{m_n\}$ of the original sequence such that all zeroes $\{z_{j,m_n}\}$ satisfy the scaling condition:
\[
\lim_{n\to\infty}\xi_{j,m_n} =\xi_j,  \qquad \xi_j \in \C^+,
\]
for every $j=1,2,\ldots,  N$, where $N=\sum_{l\ge 1}c_l\in [1,\infty]$.  Moreover,  if $N=\infty$, then $\lim_{j\to\infty}|\xi_j|=\infty$. For $z \in \D$ and $\xi = in(1-z)$, we have
\[
\frac 1n \frac{1-|z|^2}{|e^{it/n}-z|^2}
=\frac{n(1-|z|^2)}{|in(e^{it/n}-1)+in(1-z)|^2} 
=\frac{2\Im\xi - |\xi|^2/n}{|in(e^{it/n}-1)+\xi|^2}.
\]
That gives
\begin{equation}\label{sd_8}
\lim_{n\to\infty}\frac {1}{m_n} \frac{1-|z_{j,m_n}|^2}{|e^{it/m_n}-z_{j,m_n}|^2}=\frac{2\Im \xi_{j}}{|t-\xi_{j}|^2},
\end{equation}
and the convergence is uniform on compact subsets of $\R$. Since all terms in \eqref{sd_7} are non-negative, 
we can define $U$ as 
\[
U(t)=\sum_{j = 1}^{N}\frac{2\Im \xi_j}{|t-\xi_j|^2} \,.
\]
Lemma \ref{l8} and \eqref{sd_8} guarantee that
\[
\sum_{j = 1}^{l}\frac{2\Im \xi_j}{|t-\xi_j|^2} \le h'(t)=1
\]
for every $l\in \mathbb{N}$. Thus,
\[
U(t)
\le h'(t) = 1, \qquad t \in \R.
\]
Substituting $t =0$, we see that $\{\xi_j\}$ satisfies Blaschke condition in $\C^+$.  Consider the Blaschke product with zeroes at $\{\xi_j\}$, i.e.,
\[
B(\xi) =  \prod_{j=1}^N \left(e^{i\alpha_j}\frac{\xi-\xi_j}{\xi-\bar\xi_j}\right), \qquad \xi \in \C^+,
\]
where $\alpha_j$ are chosen such that 
\[
e^{i\alpha_j}\frac{i-\xi_j}{i-\bar\xi_j}>0\quad {\rm if}\quad \xi_j\neq i \quad {\rm and}\quad \alpha_j=0\quad {\rm otherwise}.
\]
We will show that $B(\xi)=e^{i\beta_1 \xi+\beta_2}$ with some $\beta_1,\beta_2\in \R$ thus getting the contradiction.
To this end,  write
\[
h'_{m_n}(t) = \frac{1}{m_n}\sum_{k=1}^{m_n}\frac{1-|z_{k,m_n}|^2}{|e^{it/m_n}-z_{k,m_n}|^2}=\Psi_{1,m_n,L}(t)+\Psi_{2,m_n,L}(t),
\]
where we define
\[
\Psi_{1,m_n,L}(t)= \frac{1}{m_n}\sum_{k \in \Omega}\frac{1-|z_{k,m_n}|^2}{|e^{it/m_n}-z_{k,m_n}|^2}\,, \quad \Omega = \{k: \; |\xi_{k}|< L,\, |\xi_{k,m_n}-\xi_k|<0.1\},
\]
and $\Psi_{2,m_n,L} = h'_{m_n} - \Psi_{1,m_n,L}$. We know from the previous lemma that $\lim_{n \to\infty} h_n' = 1$ uniformly over compacts in $\R$. When $L$ is fixed and $n\to\infty$, we have 
$$
U_L(t) \ddd \sum_{j: |\xi_j|<L}\frac{2\Im \xi_j}{|t-\xi_j|^2} =  \lim_{n\to\infty} \Psi_{1, m_n,L}(t)\le h'(t) = 1,  \qquad t \in \R.
$$
Moreover, from  \eqref{sd_22} we get $|\Psi'_{2,m_n,L}(t)|\lesssim L^{-1}\Psi_{2,m_n,L}(t)$ uniformly with respect to $t\in [-L/2,L/2]$. Since $\Psi_{1,m_n,L}+\Psi_{2, m_n,L}=h_{m_n}'\to 1$ uniformly over compacts, we have $|\Psi'_{2, m_n,L}(t)|\lesssim L^{-1}$ for $t\in [-L/2,L/2]$ if $n$ is large enough. Clearly, $\{\Psi_{2,m_n,L}\}$ converges uniformly on $[-L/2,L/2]$ as a difference of two uniformly convergent sequences. Therefore, if $\Psi_{2,L}$ denotes its limit, then 
\[
\|\Psi_{2,L}\|_{L^\infty[-L/2,L/2]}\le 1, \quad \|\Psi_{2,L}\|_{{\rm Lip}[-L/2,L/2]}\lesssim L^{-1}\,,
\]
where
\[
\|f\|_{{\rm Lip}[a,b]}\ddd\sup_{x,y\in [a,b],x\neq y}\frac{|f(x)-f(y)|}{|x-y|}\,.
\]
Notice that $U_L+\Psi_{2,L}=1$. Therefore,  $\|U_L\|_{{\rm Lip}[-L/2,L/2]}\lesssim L^{-1}$, that is
\[
|U_L(\xi_2)-U_L(\xi_1)|\lesssim |\xi_1-\xi_2|/L
\]
for $\xi_1,\xi_2\in [-L/2,L/2]$. Taking the limit as $L\to\infty$
and recalling that $U=\lim_{L\to\infty} U_L$, we see that  $U$ is constant on $\R$.
A direct calculation shows that $(\arg B)'=cU$ on $\R$ for some positive constant~$c$. Thus, we have $\arg B(t)=\beta_1t+\beta_2$, $t\in \R$, $\beta_1\ge 0$. The function $Be^{-i(\beta_1z+\beta_2)}$ is unimodular on $\R$ and has zero argument there, so it is equal to $1$ on $\R$ and, by uniqueness of holomorphic functions,  $B(z)=e^{i(\beta_1z+\beta_2)}$, $z\in \C^+$, yielding a contradiction. \qed

\medskip

\begin{Lem}\label{sd_99}
Assume that smooth functions $h_n$ defined on $(-\pi n,\pi n)$ have derivatives given by \eqref{sd_7}
and the sequence $\{h_n\}$ converges almost everywhere to some nondecreasing function $h$ defined on $\R$.  If 
\begin{equation}\label{sd_15}
\lim\limits_{n\to\infty}\left(n\min\limits_{1\le j \le n}|1-z_{j,n}|\right)=\infty\,,
\end{equation}
then $h=c_1t+c_2$ and $\{h_n\}$ converges uniformly over compacts in $\R$. 
\end{Lem}
\beginpf For arbitrary $b>0$, we have
\[
\int_{-b}^b h_n'(t)dt=h_n(b)-h_n(-b)
\]
and, since $\lim_{n\to \infty}h_n(t)= h(t)$ a.e. and $h_n'\ge 0$, one gets
\[
\sup_{n}\int_{-b}^b h_n'(t)dt<\infty\,.
\]
Moreover, condition \eqref{sd_15} and an estimate \eqref{sd_22} give
\[
\lim\limits_{n\to\infty}\int_{-b}^b |h_n''(t)|dt=0\,.
\]
From the relation 
$$
h_n'(t)(t_2-t_1)=h_n(t_2)-h_n(t_1)+\int_{t_1}^{t_2}\int_{\tau}^t h_n''(\tau_1)d\tau_1d\tau\,,
$$
we obtain
\[
\lim_{n\to \infty}h_n'(t)= \frac{h(t_2)-h(t_1)}{t_2-t_1}\,.
\]
In particular, the right hand side does not depend on $t_1$ and $t_2$. This implies that $h$ is a linear function, i.e., $h=c_1t+c_2$. Lemma \ref{sd_111} gives uniform convergence. \qed

\medskip

\section{Proof of Theorem \ref{t2}}\label{s5}
The proof of Theorem \ref{t2} is based on careful study of the arguments of orthogonal polynomials $\phi_n$ and Schur functions $f_n$. We proceed as follows. Theorem~\ref{t1}, Lemma~\ref{l7}, and Lemma~\ref{l12}  show that these arguments, after rescaling and taking a limit, satisfy equation \eqref{eq4} below. Since the derivative  of the argument of a polynomial with zeroes in $\D$ is a finite sum of Poisson kernels (see formula \eqref{sd_a21}), equation \eqref{eq4} allows us to recover local asymptotics of all objects in Theorem \ref{t2} and prove that assertions $(a)$--$(d)$ are equivalent to an identity $d=0$ in \eqref{eq4}. 

\medskip

In this section, we always assume $\mu \in \szc$. We start with several auxiliary results. Recall that $\{f_n\}$ denotes the family of Schur functions for a measure $\mu$ and that, given an arc  $I\subset \T$ with the center at $\xi_I\in \T$, we let $z_I=(1 - |I|)\xi_I$. 
\begin{Lem}\label{l14}
Let $I \subset \T$ be an arc  and $|I| \le 1/4$. Then,
$$
|\{\xi \in I: |f_n(\xi) - f_n(z_I)| > t\}|/|I| \lesssim \eta(z_I)/t,
$$
where the function $\eta$ is defined in \eqref{sd_a1}. 
\end{Lem}
\beginpf By  Lemma \ref{l7}, we have $\|f_n\|^{*}_{\eta} \le C$ for some constant $C$. Since 
$$
\frac{1 -|z_I|^2}{|1 - \bar \xi z_I|^2} \gtrsim \frac{1 - (1-|I|)^2}{|I|^2} \gtrsim \frac{1}{|I|}, \qquad \xi \in I,
$$ one has
$$
 \frac{1}{|I|}\int_{I}|f_n(\xi) - f_n(z_I)|\,dm(\xi) \lesssim  \P(|f_n - f_n(z_I)|, z_I) \le C\eta(z_I).
$$
It remains to use Chebyshev inequality.\qed

\medskip

Given $\xi \in \T$, $\rho\in (0,1)$, $\delta\in (0,1)$ and $\alpha$, $\beta$: $0<\alpha<\beta$, set 
\begin{equation}\label{sd_a66}
\Upsilon_{\delta,\rho,\alpha,\beta}(\xi)=\{z\in S^*_\rho(\xi), \; \alpha\delta<|z-\xi|<\beta\delta\},
\end{equation}\smallskip
where, as before, $S^*_\rho(\xi)$ is the convex hull of $\rho\D$ and point $\xi$. For a complex-valued function $h$ defined on a domain $\Omega \subset \C$, we introduce its oscillation as
\[
{\rm osc}_{\Omega}(h)=\sup_{z_1,z_2\in \Omega}|h(z_2)-h(z_1)|.
\]
In the next lemma, we show that Schur family $\{f_k\}$ has  small oscillation near the boundary of $\D$ uniformly in $k \ge 0$.
\begin{Lem}\label{l11}
Suppose $\xi \in \T$ is such that $\lim_{r\to 1}\K(\mu,r\xi) = 0$. Then, for every $\rho$,$\alpha$,$\beta$ and $\{\delta_n\}$ such that $\lim_{n\to\infty}\delta_n=0$, we have
$$
\lim_{n \to +\infty} \sup_k{\rm osc}_{\Upsilon_n} (f_k)= 0, \qquad \Upsilon_n = \Upsilon_{\delta_n,\rho,\alpha,\beta}(\xi).
$$ 
\end{Lem}
\beginpf Take an arc $I_n \subset \T$ centered at $\xi$ so that $|I_n|=c_n\delta_n$ for some $c_n > 0$ such that $c_n\to\infty$, $c_n\delta_n \to 0$, and $c_n\sqrt{\rho_n} \to 0$, where $\rho_n=\eta(z_n)$, $z_n=\xi (1-|I_n|)$. For example, one can take $c_n = 1/(\sqrt{\delta_n} + \sqrt[4]{\widetilde\rho_n})$,  $\widetilde \rho_n = \sup_{r\ge 1-\sqrt{\delta_n}}\eta(r\xi)$. By Lemma~\ref{l14}, we have
\begin{equation}\label{s_1}
|\{\xi \in I_n: |f_k(\xi) - f_k(z_n)| > t\}|/|I_n| \lesssim \eta(z_n)/t, \quad t>0.
\end{equation}
For every $g$ that satisfies $\|g\|_{L^\infty(\T)}\le 2$ and every $z\in\Upsilon_{n}$, we have
\[
\P(g, z)=\P(\chi_{I_n}g, z)+\P(\chi_{\T\setminus I_n}g, z) = \P(\chi_{I_n}g, z) + o(1),\quad n\to\infty,
\]
since $\lim_{n\to\infty}c_n=\infty$, and this bound holds uniformly in $g$ and $z$. Thus, having defined $\widetilde f_k=(f_k-f_k(z_n))\chi_{I_n}$, we get
\begin{equation}\label{sd_a5}
f_k(z)-f_k(z_n)=\P(f_k- f_k(z_n), z)=\P(\widetilde f_k, z)+o(1).
\end{equation}
Recall that $\lim_{n \to \infty}c_n\sqrt{\rho_n} = 0$. Thus,
\begin{align*}
\left|\P(\widetilde f_k, z)\right| &\le \P(\chi_{|\widetilde f_k|<\sqrt\rho_n}|\widetilde f_k|, z) + \P(\chi_{|\widetilde f_k|\ge\sqrt\rho_n}|\widetilde f_k|, z)\,.
\end{align*}
The first term is bounded by $\sqrt{\rho_n}$.
Consider the second one. Since $z\in\Upsilon_{n}$, we can estimate the Poisson kernel by $C\delta_n^{-1}$, bound $|\widetilde f_k|$ by $2$, and apply \eqref{s_1} to write
\begin{align*}
\P(\chi_{|\widetilde f_k|\ge\sqrt\rho_n}|\widetilde f_k|, z)
\lesssim  \frac{|I_n|\rho_n}{\sqrt\rho_n\delta_n}=c_n\sqrt{\rho_n} \to 0.
\end{align*}
From \eqref{sd_a5}, we get
\[
\lim_{n\to\infty}\sup_k\sup_{z\in \Upsilon_{n}}|f_k(z)-f_k(z_n)|=0\,.
\]
Since $|f_k(\xi_1)-f_k(\xi_2)|\le  |f_k(\xi_1)-f_k(z_n)|
+|f_k(\xi_2)-f_k(z_n)|$, we get the statement of the lemma. \qed
 
\medskip

Denote the argument of $\phi_n^*$ on $\T$ by $\zeta_n$. Since $\phi_n^*$ has no zeroes in $\overline{\mathbb{D}}$, $\zeta_n=\Im \log \phi_n^*(e^{it})$ is a continuous function and it coincides with the harmonic conjugate of $\log |\phi_n^*(e^{it})|$ since $\phi_n^*(0)$ is real. Moreover, $\zeta_n(e^{it}) = (nt - \gamma_n(t))/2$ where $\gamma_n$ denotes an argument of the Blaschke product $b_n = \phi_n/\phi_n^*$. As was discussed previously, $\gamma_n(t)$ is increasing in $t\in [-\pi,\pi)$, see \eqref{sd_a21}.

\begin{Lem}\label{l17}
The function $\phi_n^*(1-zb_nf_n)$ is outer in $\D$. For almost every $t \in (-\pi, \pi)$, the harmonic conjugate of the function $\log|\phi_n^*(1-\xi b_nf_n)|^2$, $\xi\in \T$, at point $e^{it}$ is given by
\begin{equation}\label{eq5}
v_n(t) = nt - \gamma_n(t) + 2\arctan\left(\frac{|f_n(e^{it})|\sin(\gamma_n(t) + t +  \kappa_n(t)) }{1 + |f_n(e^{it})| \cos(\gamma_n(t)+ t + \kappa_n(t))}\right)\,.
\end{equation} 
In this formula, the function $\kappa_n(t)$ is uniquely defined by conditions: $\kappa_n(t)\in [-\pi,\pi)$  and $e^{i\kappa_n(t)}=-f_n(e^{it})/|f_n(e^{it})|$ in the case when $f_n$ is not identically zero.  If $f_n=0$ identically, then the third term in \eqref{eq5} can be dropped.
\end{Lem}
\beginpf Since the polynomial $\phi_n^*$ has no zeroes in $\D$ it is an outer function. We also have $\Re (1-zb_nf_n) \ge 0$ in $\D$, hence $1-zb_nf_n$ is an outer function as well, see Corollary 4.8 on page 74 in \cite{Garnett}. The harmonic conjugate of $\log|\phi_n^*(1-zb_nf_n)|$ is the sum of harmonic conjugates of $\log|\phi_n^*|$ and of $\log|1-zb_nf_n|$. The harmonic conjugate of $\log|\phi_n^*|$ is $\zeta_n$. The harmonic conjugate of $g = \log |1-\xi b_nf_n|$ is equal to $\Im \log(1-\xi b_nf_n)$ which is the boundary value of the argument of function $1-zb_nf_n$. The latter function has positive real part and its absolute value is bounded by $2$. Therefore, $\widetilde g$ is well defined a.e. on $\T$  and $\widetilde g\in [-\pi/2,\pi/2]$.
As we have seen in \eqref{eq08}, for $\mu\in \sz$ we have 
\[
\int_\T \log (1-|f_n|^2)\,dm < +\infty,
\]
in particular, $|f_n|<1$ almost everywhere on $\T$ for each $n$. Suppose that $\xi = e^{it}$ is such that $f_n(\xi)$, the boundary value of $f_n$, satisfies $0<|f_n(\xi)|<1$. We know that this holds for almost every $\xi=e^{it}\in \T$ (if $|f_n|=0$ on a set of positive Lebesgue measure, then $f_n=0$ identically and the lemma holds trivially). Take $\kappa_n(t)\in [-\pi,\pi)$ such that  $-f_n(e^{it})/|f_n(e^{it})|=e^{i\kappa_n(t)}$. Then, we have
$$
\widetilde g(\xi) = \arctan\left(\frac{|b_n(\xi)f_{n}(\xi)|\sin(\gamma_n(t)+t+\kappa_{n}(t))}{1+|b_n(\xi)f_{n}(\xi)|\cos(\gamma_n(t)+t+\kappa_{n}(t))}\right),
$$
due to the formula
\begin{equation}\label{eq3}
\frac{1+ae^{i\psi}}{|1+ae^{i\psi}|} = \exp\left(i\arctan\left(\frac{a\sin\psi}{1+a\cos\psi}\right)\right), 
\end{equation}
for $a\in [0,1]$, $\psi \in \R$: $1+a\cos\psi \neq 0$, when we notice that $\Re(1+\xi b_nf_n) \neq 0$ almost everywhere on $\T$. Since $|b_n(\xi)|=1$, the lemma is proved.\qed

\medskip Now, we can control the oscillation of $\upsilon_n$.

\begin{Lem}\label{l15}
Let $v_n$ be defined by \eqref{eq5} and let $I \subset \T$ be an arc with center at $\xi_0 \in \T$. Then, there exist numbers $c_{I,n}$ such that
$$
|\{\xi \in \T: |v_n(\xi) - c_{I,n}| > t\}|/|I| \lesssim  t^{-1}\sum_{j \ge 0}2^{-j}\eta(z_j),
$$
where the function $\eta$ is defined in \eqref{sd_a1}  and $\{z_j\}$ is the set of points constructed in Lemma \ref{l12}.
\end{Lem}
\beginpf Since $v_n$ is the harmonic conjugate of  $u_n = \log|\phi_n^*(1 - zb_nf_{n})|^2$,  we obtain
$$
|\{\xi \in \T: |v_n(\xi) - c_{I,n}| > t\}|/|I| \lesssim  t^{-1}\|u_n\|_{ \eta}^{*} \sum\nolimits_{j\ge 0}2^{-j}\eta(z_j)
$$
from Lemma~\ref{l12}.
It remains to note that $\{\|u_n\|_{\eta}^{*}\}$ is uniformly bounded due to Lemma~\ref{l7}. \qed

\medskip
\noindent We recall that Christoffel-Darboux kernel is defined by
\[
k_{\xi,\mu,n}(z)=\sum_{j=0}^{n-1}\phi_j(z)\overline{\phi_j(\xi)}\,,
\]
where $\{\phi_j\}$ are polynomials orthonormal with respect to measure $\mu$.

\begin{Lem}\label{l1}
If $\xi = e^{it}$ and $t\in \R$, then $\|k_{\xi, \mu, n}\|^{2}_{L^2(\mu)} = |\phi_n^*(\xi)|^2\gamma'_n(t)$.
\end{Lem}
\beginpf For $\xi \in \T$ and $z \neq \xi$, we have (see  \cite{ger}, Section 1)
$$
k_{\xi,\mu,n}(z) = \frac{\phi_n^*(z)\ov{\phi_n^*(\xi)} - \phi_n(z)\ov{\phi_n(\xi)}}{1 - z \bar \xi} = \phi_n^*(z)\ov{\phi_n^*(\xi)}\frac{1 - b_n(z)\ov{b_n(\xi)}}{1 - z \bar \xi}.
$$  
Noting that $b_n(e^{is}) = e^{i\gamma_n(s)}$ for $s \in [-\pi, \pi)$, we get 
\begin{align*}
\|k_{\xi,\mu,n}\|_{L^2(\mu)}^{2} 
&=k_{\xi,\mu,n}(\xi)= \lim_{z \to \xi} k_{\xi,\mu,n}(z), \\
&= |\phi_n^*(\xi)|^2\lim_{s \to t}\frac{1 - e^{i(\gamma_n(s) - \gamma_n(t))}}{1 - e^{i(s-t)}}, \\
&= |\phi_n^*(\xi)|^2\gamma'_n(t). 
\end{align*}
The lemma follows. \qed

\medskip

\begin{Lem}\label{lp1}
Assume that $\lim_{n\to\infty}|\phi_n^*(\xi)|^{-2} = |D_{\mu}(\xi)|^{2}$ for almost every $\xi \in \T$. Let $r_n=1-1/n$ for $n \ge 1$. Then,  $\lim_{n\to\infty} f_n(r_n \xi) = 0$ for almost every $\xi \in \T$.
\end{Lem}
\beginpf  We claim that for each $\delta \in (0,1)$, there exists a subset $G_{\delta}(\mu) \subset \T$ with the properties:
\begin{align}
&m(G_{\delta}(\mu)) \ge 1- \delta,  \label{A}\\
&\mbox{each point of } G_{\delta}(\mu) \mbox{ is a Lebesgue point,}\label{C}\\
&\lim\limits_{\eps\to 0}\sup\{\K(\mu,z),\; z\in S_{\rho}^*(\xi),\, |z-\xi|<\eps\}=0 \mbox{ for } \xi \in G_\delta(\mu) \label{D}\, \mbox{and } \rho\in (0,1),\\
&\lim\limits_{n \to\infty} n^{-1}\cdot\|k_{\xi, \mu, n}\|^{2}_{L^2(\mu)} = |D_{\mu}(\xi)|^{-2}\mbox{ for } \xi \in G_\delta(\mu), \label{sd_a11}\\
&\lim\limits_{n \to \infty}\gamma'_n(t)/n = 1 \mbox{ uniformly with respect to }t: \; e^{it} \in G_{\delta}(\mu).\label{sd_2}
\end{align}
Indeed, for every $\mu \in \szc$ we have \eqref{D} almost everywhere on $\T$ since the Poisson kernel is an approximate identity. Theorem 1 in \cite{MNT91} says that the limit relation in \eqref{sd_a11} holds almost everywhere on $\T$. By Lemma \ref{l1}, this implies that the limit relation in \eqref{sd_2} holds almost everywhere on $[-\pi,\pi]$. So, there is a set $E$ of full Lebesgue measure on $\T$ such that limit relations in \eqref{sd_a11}, \eqref{D}, \eqref{sd_2} holds for $\xi = e^{it}$ in $E$. Using Egorov's theorem, one can find a subset $\widetilde E$ of $E$ of length $2\pi(1-\delta)$ such that the limit relation in \eqref{sd_2} is uniform with respect to $t \in [-\pi, \pi]$ provided $e^{it} \in \widetilde E$. Then, we can denote by $G_{\delta}(\mu)$ the set of the Lebesgue points of $\widetilde E$.

\medskip

Now, it suffices to prove that for every fixed $\delta>0$ we have $\lim_{n\to\infty}f_n(r_n \xi) = 0$ for every $\xi \in G_{\delta}(\mu)$. 
Without loss of generality, we can assume that $\xi = 1$. Consider any convergent subsequence $\{f_{n_k}(r_{n_k})\}$ and let $\lim_{k \to \infty} |f_{n_k}(r_{n_k})| = d$. Let us show that $d = 0$. During this proof, we will several times choose subsequences of $\{f_{n_k}(r_{n_k})\}$. To simplify notation, we will assume that sequences under consideration converge without extracting subsequences. In particular, we let $\lim_{n \to \infty} |f_n(r_n)| = d$. By Lemma \ref{l11}, we have $\lim_{n \to \infty} |f_n(1-a/n)| = d$ for every $a>0$. Let, as before, $\gamma_n: [-\pi, \pi) \to \R$ be a continuous branch of the argument of the function $b_n(e^{it})$, where $b_n = \phi_n/\phi^*_n$. Denote 
$I_n(a) = (-a/n, a/n)$ for all $n \ge 1$ and a constant $a \ge 10$. Consider $n \ge a/\pi$. It follows from Lemma~\ref{l17}, Lemma~\ref{l15}, and condition \eqref{D} that there are sets $E_{n}(a) \subset I_n(a)$ and numbers $c_n$ such that functions 
\begin{equation}\label{eq2}
v_n(t) = nt - \gamma_n(t) + 2\arctan\left(\frac{|f_n(e^{it})|\sin(\gamma_n(t) + t + \kappa_n(t)) }{1 + |f_n(e^{it})| \cos(\gamma_n(t)+t+\kappa_n(t))}\right)
\end{equation}
satisfy the following relations:
\begin{itemize}
\item[$(a)$] $|v_n(t) - c_{n}| \le \eps_{n}$ for all $t \in E_{n}(a)$, 
\item[$(b)$] $|E_{n}(a)| \ge (1-\eps_{n})|I_{n}(a)|$, 
\end{itemize}  
for some positive sequence  $\{\eps_n\}_{n \ge 1}$ converging to zero. Next, we renormalize \eqref{eq2} as follows. For each $n$, take $\pi_n\in \{2\pi \Z\}$ such that $|c_n-\pi_n|\le \pi$ so 
\[
|(v_n(t)-\pi_n)-(c_n-\pi_n)|\le \eps_n
\]
for all $t\in E_n(a)$. We denote $\widehat c_n=c_n-\pi_n$, $\widehat v_n=v_n-\pi_n$ and $\widehat\gamma_n=\gamma_n+\pi_n$. Now, \eqref{eq2} can be rewritten as
\begin{equation}\label{eq2n}
\widehat v_n(t) = nt - \widehat \gamma_n(t) + 2\arctan\left(\frac{|f_n(e^{it})|\sin(\widehat\gamma_n(t) + t + \kappa_n(t)) }{1 + |f_n(e^{it})| \cos(\widehat \gamma_n(t)+t+\kappa_n(t))}\right)
\end{equation}
and the following relations hold:
\begin{itemize}
\item[$(a')$] $|\widehat c_n|\le \pi$ and $|\widehat v_n(t) - \widehat c_{n}| \le \eps_{n}$ for all $t \in E_{n}(a)$.
\end{itemize}
Since $\widehat \gamma_n$ is increasing on $(-\pi,\pi)$, relations \eqref{eq2n} and $(a')$ imply that there is a constant $c(a)$ depending only on $a$, such that
$|\widehat\gamma_n(t)| \le c(a)$, $t \in \mathrm{co}\,(E_n(a))$, where $\mathrm{co}\,(E_n(a))$ is the convex hull of the set $E_{n}(a) \subset I_n(a)$. Note that $\mathrm{co}\,(E_n(a))$ contains $I_n(a/2)$ for $n$ such that $|\eps_n| \le 1/2$. Hence, for large enough  $n$, the functions 
$$
h_n: s \mapsto \widehat\gamma_n(s/n)
$$ 
are correctly defined on $[-a/2,a/2]$, increasing, and uniformly bounded by $c(a)$. Therefore, by Helly's selection theorem, one can choose a subsequence of $\{h_n\}$ that converges pointwise on $[-a/2,a/2]$ to a non-decreasing function $h$. We again will assume that the whole sequence converges to $h$. One can also assume that functions $\widehat v_n(s/n)$, $|f_n(e^{is/n})|$ on $[-a/2,a/2]$ converge in measure to constants $c \in [-\pi, \pi]$, $d$, respectively. Indeed, for $\widehat v_n(s/n)$ this follows from assertion $(a')$, while for $|f_n(e^{is/n})|$ -- from Lemma \ref{l14}. If $d \neq 0$, Lemma \ref{l14} implies also the converge of $\kappa_n(s/n)$ in measure on $[-a/2,a/2]$ to a constant $\kappa \in [-\pi, \pi]$.  Choosing, if needed, a subsequence, one can assume  (see \cite{folland}, Theorem 2.30) that the convergence of $\widehat v_n(s/n)$, $|f_n(e^{is/n})|$ and $\kappa_n(s/n)$ is pointwise on a subset $E \subset [-a/2,a/2]$ of full Lebesgue measure. Since $\xi=1$ is the Lebesgue point of the set $G_{\delta}(\mu)$ and $h_n'(s) = \gamma'_n(s/n)/n$, we use \eqref{sd_2} to get 
\begin{equation}\label{sd_3}
h(s_2) - h(s_1) = \lim_{n \to +\infty} (h_n(s_2) - h_n(s_1)) = \lim_{n \to +\infty} \int_{s_1}^{s_2}h_n'(s)\,ds \ge s_2 - s_1   
\end{equation}
for every $s_1 \le s_2$ in $E$. We consider two cases now.

\medskip

\noindent {\bf Case 1.} If $d \in [0,1)$, then relation \eqref{eq2n} implies
\begin{equation}\label{eq4}
c = s - h(s) + 2\arctan\left(\frac{d\cdot \sin (h(s)+\kappa) }{1 + d\cdot \cos (h(s)+\kappa)}\right), \qquad s \in E,
\end{equation}
where we set $\kappa = 0$ if $d = 0$. The derivative
\[
\partial_h\left( h-2\arctan \left( \frac{d\sin (h+\kappa)}{1+d\cos(h+\kappa)}  \right)\right)=\frac{1-d^2}{1+2d\cos(h+\kappa)+d^2}
\]
is within $[(1-d)/(1+d),(1+d)/(1-d)]$ so application of the inverse function theorem shows that \eqref{eq4} defines a smooth increasing function  on $[-a/2,a/2]$. Since $h$ is nondecreasing, we see that \eqref{eq4} holds for all $s \in [-a/2,a/2]$. Moreover, \eqref{sd_3} gives  $h'(s) \ge 1$ for such $s$. Differentiating \eqref{eq4}, we obtain
$$
h'(s) = \frac{1 + 2d\cos (h(s)+\kappa) + d^2}{1 - d^2}, \qquad s \in [-a/2,a/2].
$$ 
Since $h'\ge 1$ and a parameter $a$ is large enough, there is 
 $s^\ast \in [-a/2,a/2]$ so that $\cos (h(s^\ast)+\kappa) = -1$. Thus, $(1-d)/(1+d) \ge 1$ which implies $d = 0$ and we are done.

\medskip

\noindent {\bf Case 2.} Let $d = 1$ and rewrite \eqref{eq2n} as
\begin{align}
\widehat v_n(s/n) 
&= s - \widehat \gamma_n(s/n) + 2\arctan\left(\frac{|f_n(e^{is/n})|\sin(\widehat \gamma_n(s/n)+s/n+\kappa_n(s/n))}{1 + |f_n(e^{is/n})| \cos(\widehat \gamma_n(s/n)+s/n+\kappa_n(s/n))}\right). \label{eq2nn}
\end{align}
Taking the limit requires some care in this case. We have  
\[
\lim_{n\to\infty}\Bigr(1-e^{is/n}f_n(e^{is/n})b_n(e^{is/n})\Bigl)= 1+e^{i(\kappa+h(s))}
\]
for almost every $s\in [-a/2,a/2]$. Let $\widetilde E$ be a subset of $E$ on which $1+e^{i(\kappa+h(s))}\neq 0$. 
If the function $H$ is defined by the formula
\[
H(\alpha)=(\alpha-2\pi j)/2\qquad {\rm if}\; \alpha\in (2\pi j-\pi,2\pi j+\pi), \quad j\in \Z, 
\]
then an identity
\begin{equation}\label{sad_shame1}
\arctan \left(\frac{\sin\alpha}{1+\cos\alpha}\right)=H(\alpha), \quad \alpha: \;\cos\alpha \neq -1,
\end{equation}
is immediate. Given \eqref{sad_shame1},  take a limit in \eqref{eq2nn} for every $s\in \widetilde E$ to get
\[
c = s - h(s) + 2\arctan\left(\frac{\sin (h(s)+\kappa) }{1 + \cos(h(s)+\kappa)}\right)=s-h(s)+2H(h(s)+\kappa)\,.
\]
Thus, if $s_1\neq s_2$ and $s_1,s_2\in \widetilde E$, then $s_2-s_1\in \pi\Z$ and so $\widetilde E$ is either finite or empty. This implies that $e^{i(\kappa +h)}=-1$ almost everywhere on $[-a/2,a/2]$ and $h$ is a nondecreasing step function. That, however, contradicts \eqref{sd_3} and we get $d\neq 1$ under assumptions of the lemma. \qed
 
\medskip

We recall that the zeroes of $\phi_n$ were denoted by $\{z_{j,n}\}$ and they are all inside $\D$.
 \begin{Lem} \label{sd_33} Suppose there is $a>0$ such that
 \[
 \lim_{n\to\infty}f_n(r_{a,n}\xi)=0
 \]
for almost every $\xi\in \T$. Then, 
\begin{equation}\notag
 \lim_{n\to\infty}|\phi_n^*(\xi)|^2= |D(\xi)|^{-2}
 \end{equation}
  and 
 \begin{equation}\label{sd_10}
\liminf\limits_{n\to\infty} \Bigl(n\min_{1\le j\le n}|\xi-z_{j,n}|\Bigr)=+\infty
 \end{equation}
 for almost every $\xi\in \T$.
 \end{Lem}
\beginpf Notice that we have
\begin{equation}\label{sd_aew}
\lim_{n\to\infty}\sup_{z\in \Upsilon_{n^{-1},\rho,\alpha,\beta}(\xi)}|f_n(z)|=0
\end{equation}
for all $\rho,\alpha,\beta$ and almost every $\xi\in \T$ by Lemma \ref{l11}. Moreover, for almost every $\xi\in \T$, we have
\begin{equation}\label{sd_a121}
\lim_{\epsilon\to 0}\sup_{z\in S_{\rho}^*(\xi),|z-\xi|<\epsilon}\K(\mu,z)=0
\end{equation}
for any  $\rho\in (0,1)$, and 
\begin{equation}
\lim_{n \to\infty} n^{-1}\cdot\|k_{\xi, \mu, n}\|^{2}_{L^2(\mu)} = |D_{\mu}(\xi)|^{-2}\,.\label{sd_a114}
\end{equation}
Without loss of generality, we assume that $\xi=1$ is a point at which all these conditions are satisfied and let $I_n(a)= [-a/n,a/n]$. Like in the proof of previous lemma, we have
\begin{equation}\label{eq21}
\widehat v_n(t) = nt - \widehat \gamma_n(t) + 2\arctan\left(\frac{|f_n(e^{it})|\sin(\widehat\gamma_n(t) + t + \kappa_n(t)) }{1 + |f_n(e^{it})| \cos(\widehat\gamma_n(t)+t+\kappa_n(t))}\right)
\end{equation}
and there is a set $E_n(a)\subseteq I_n(a)$ and a sequence $ \widehat c_n\in [-\pi,\pi)$ such that 
\begin{itemize}
\item[$(a')$] $|\widehat v_n(t) - \widehat c_{n}| \le \eps_{n}$ for all $t \in E_{n}(a)$, 
\item[$(b)$] $|E_{n}(a)| \ge (1-\eps_{n})|I_{n}(a)|$, 
\end{itemize}  
for some positive sequence  $\{\eps_n\}_{n \ge 1}$ converging to zero. Moreover, we can use \eqref{sd_aew}, Lemma \ref{l14}, and \eqref{sd_a121} to choose $E_n(a)$  such that an additional condition
\begin{itemize}
\item[$(c)$] $|f_n(e^{it})|\le \eps_n$ for all $t\in E_n(a)$
\end{itemize}
is satisfied.
Rescale $t\in I_n(a)$ as $t=s/n$, let $h_n(s)=\widehat \gamma_n(s/n)$ as in the previous proof, and write
$$
\widehat v_n(s/n) = s - h_n(s) + 2\arctan\left(\frac{|f_n(e^{is/n})|\sin(h_n(s) + s/n + \kappa_n(s/n)) }{1 + |f_n(e^{is/n})| \cos(h_n(s)+\tau/n+\kappa_n(s/n))}\right).
$$
Take a  limit in measure on $[-a/2,a/2]$ in the above equation through some subsequence $\{n_j\}$, it exists thanks to $(a')-(c)$. It follows from $(c)$ that the sequence $\{h_{n_j}\}$ converges to $s+c$ in measure, where $c$ can depend on the choice of subsequence $\{n_j\}$. From each functional sequence converging in measure, we can choose a subsequence converging almost everywhere. We denote it by the same $\{n_j\}$. Since each function $h_n$ is increasing,  this convergence  is in fact uniform over $[-a/2,a/2]$ due to Lemma~\ref{sd_111}. The parameter $a$ was arbitrary so we can take an unbounded positive sequence $\{a_l\}$ and choose the subsequence of $\{h_n\}$ which  converges to a linear function uniformly on all compacts in $\R$. We again denote it by $\{h_{n_j}\}$. Thus, in view of \eqref{sd_a21}, we can apply Lemma~\ref{l8}  to show that $\{h_{n_j}'\}$ converges to $1$ uniformly over compacts and this convergence holds at point $s=0$, in particular. Arguing by contradiction, we can prove that in fact $\lim_{n\to\infty}h'_n(0)=1$ through the whole sequence. By \eqref{sd_a114} and  Lemma \ref{l1}, we get an implication:
\begin{equation}\label{sd_im}
\lim_{n\to\infty}\frac{\gamma_n'(0)}{n} = 1 \!\implies\! \lim_{n\to\infty}|\phi^*_n(1)|^2 =|D_\mu(1)|^{-2}.
\end{equation}
The property  \eqref{sd_10} of zeroes follows from Lemma~\ref{sd_12}. Indeed, $\lim_{n\to\infty}h_n'= 1$ uniformly over compacts in $\R$ so every subsequential limit of $\{h_n\}$ is a linear function. So, if 
\begin{equation}\label{sd_1011}
\liminf\limits_{n\to\infty} \Bigl(n\min_{1\le j\le n}|1-z_{j,n}|\Bigr)<\infty\,,
\end{equation}
we can choose a subsequence $\{k_n\}$ over which, first, $\{h_{k_n}\}$ converges uniformly to a linear function and, secondly,
\begin{equation}\label{sd_10111}
\liminf\limits_{n\to\infty} \Bigl(k_n\min_{1\le j\le n}|1-z_{j,k_n}|\Bigr)<\infty\,.
\end{equation}
That contradicts Lemma \ref{sd_12}. \qed

\medskip

The next result shows that information about zeroes $\{z_{j,n}\}$ gives  control of pointwise asymptotics of $\{|\phi_n(\xi)|\}$ for $\xi\in \T.$

\begin{Lem}\label{sd_44} Suppose that
\begin{equation}\label{sd_13}
\lim\limits_{n\to\infty} \Bigl(n\min_{1\le j\le n}|\xi-z_{j,n}|\Bigr)=+\infty
 \end{equation}
holds for almost every $\xi\in \T$. Then, 
\begin{equation}\notag
 \lim_{n\to\infty}|\phi_n^*(\xi)|^2= |D_\mu(\xi)|^{-2}
 \end{equation}
almost everywhere on $\T$.
\end{Lem}
\beginpf We consider $\xi$ in the full measure set of points on $\T$ where \eqref{eq2n} and \eqref{sd_13} hold. Assume again without loss of generality, that $\xi=1$ and write renormalized equation \eqref{eq2n} taking $s=tn$
\begin{equation}\label{eq212}
\widehat v_n(s/n) = s - h_n(s) + 2\arctan\left(\frac{|f_n(e^{is/n})|\sin(h_n(s) + s/n + \kappa_n(s/n)) }{1 + |f_n(e^{is/n})| \cos(h_n(s)+s/n+\kappa_n(s/n))}\right)
\end{equation}
and 
\begin{itemize}
\item[$(a')$] $|\widehat v_n(t) - \widehat c_{n}| \le \eps_{n}$ for all $t \in E_{n}$, and $|\widehat c_n|\le \pi$.
\item[$(b)$] $|E_{n}(a)| \ge (1-\eps_{n})|I_{n}(a)|$, 
\end{itemize}  
for some positive sequence  $\{\eps_n\}_{n \ge 1}$ converging to zero. Therefore, since $h_n$ is increasing,
\[
\sup\{|h_n(s)|, \; n \ge 1, \; s\in [-a/2,a/2]\} <\infty,
\]
and we can apply Helly's theorem on $[-a/2,a/2]$ to find a  subsequence $\{h_{k_n}\}$ which converges to a limit $h$ almost everywhere on $[-a/2,a/2]$. The parameter $a$ is arbitrary so, going to subsequences, we can find a non-decreasing function $h$ defined on $\R$ such that a subsequence of $\{h_{n}\}$ (call it $\{h_{k_n}\}$ also) converges to $h$ almost everywhere on  $\R$.  From Lemma~\ref{sd_99}, we know that $h=c_1t+c_2$ and convergence $\lim_{n\to\infty}h_{k_n} = h$ is in fact uniform on compact subsets of $\R$. Formula \eqref{eq212} gives $c_1=1$ if we compare the variations of both sides on $[-a,a]$ when $a\to \infty$. Now, by Lemma~\ref{l8}, we get $\lim_{n\to\infty}h_{k_n}'= 1$ uniformly over compacts in $\R$. In particular, $\lim_{n\to\infty} h_{k_n}'(0) = 1$. Arguing by contradiction, we again can show that $\lim_{n\to\infty} h_{n}'(0) = 1$ over the whole sequence. By the same reasoning we used in \eqref{sd_im}, one gets the statement of the lemma.\qed

\medskip

What we proved so far implies that assertions $(a)$, $(b)$, $(c)$ of Theorem \ref{t2} are equivalent on a subset of $\T$ of full Lebesgue measure. Let us proceed with item $(d)$. The paper \cite{KH01} will be the main reference in many arguments given below.

\medskip

\begin{Lem}\label{sd_117}
Suppose $\lim_{n\to\infty}f_n=0$ almost everywhere on $\T$. Then, 
$$
\lim_{n\to\infty}|\phi_n^*|^2=|D_\mu|^{-2}
$$ 
almost everywhere on $\T$.
\end{Lem}
\beginpf This immediately follows from Khrushchev's formula 
\[
|\phi_n^*|^2|D_\mu|^{2}=\frac{1-|f_n|^2}{|1-\xi b_nf_n|^2},
\]
see identity (1.18) in \cite{KH01}. \qed

\medskip

Given $\mu$, we recall that the dual measure $\mu_{\rm dual}$ corresponds to the Schur function which is equal to $-f$. The associated orthonormal polynomials are called the polynomials of the second kind and they are denoted $\{\psi_n\}$. The Wall polynomials $\{A_n\},\{B_n\}$ are connected to orthogonal polynomials by (see formula (5.5) in \cite{KH01})
\begin{eqnarray}
\phi_{n+1}=k_{n+1}(zB_n^*-A_n^*), \quad \phi_{n+1}^*=k_{n+1}(B_n-zA_n)\,,\label{sd_u}\\
\psi_{n+1}=k_{n+1}(zB_n^*+A_n^*),\quad \psi_{n+1}^*=k_{n+1}(B_n+zA_n)\,,\label{sd_u2}
\end{eqnarray}
where $k_n$ is the leading coefficient of $\phi_n$. For $n \ge 1$, let $\widehat f_n$ be the Schur function of the probability measure $|\phi_{n}^*|^{-2}\,dm$. In fact, we have $\widehat f_n = A_{n-1}/B_{n-1}$, see formula (5.11) in \cite{KH01}. 
Define 
$$
F=\frac{1+zf}{1-zf}, \qquad \widehat F_n=\frac{1+z\widehat f_n}{1-z\widehat f_n}.
$$ 
Then, from  \eqref{sd_u}, \eqref{sd_u2} we get identities (see also formulas (5.10) and (5.11) in \cite{KH01}):
$$
\widehat F_n=  \frac{1+zA_{n-1}/B_{n-1}}{1-zA_{n-1}/B_{n-1}} = \frac{\psi_n^*}{\phi_n^*}.
$$ 
To show that $(d)$ in Theorem \ref{t2} follows from the other conditions, we proceed as follows. Since the real part of $\widehat F_n$ is nonnegative, it is an outer function in $\D$ and its behavior can be controlled by the argument. We have $\arg \widehat F_n = \arg \psi_n^* - \arg \phi_n^*$ and this identity will give $\lim_{n\to\infty}\widehat F_n(\xi)= F(\xi)$. The latter condition implies $\lim_{n\to\infty}\widehat f_n(\xi)= f(\xi)$, which yields $\lim_{n\to\infty}f_n(\xi)=0$ by lemma 4.8 in \cite{KH01}. The estimate on the nontangential maximal function will easily follow.

\begin{Lem}\label{sd_55}
Suppose $Z_n\subset \overline{\D}$ and $\lim_{n\to \infty}\sup_{z\in Z_n}|f_n(z)|= 0$, then 
\[
\lim_{n\to\infty}\sup_{z\in Z_n}|f(z)-\widehat f_n(z)|= 0.
\]
\end{Lem}
\beginpf Formula (4.19) in \cite{KH01} reads
\[
f=\frac{A_n+zB_n^*f_{n+1}}{B_n+zA_n^*f_{n+1}}\,.
\]
That yields
\[
|f - \widehat f_{n+1}| = |f - A_n/B_n|=\left|\frac{f_{n+1}z({B_n^*}/{B_n}-(A_n^*A_n)/(B_n^2))}{1+zf_{n+1}A_n^*B_n^{-1}}\right|\,.
\]
We have $|A_n^*/{B_n}|\le 1, |A_n/B_n|\le 1$ in $\D$ (see Lemma 4.5 in \cite{KH01}). Moreover, since $B_n$ does not vanish in $\overline\D$ (by the same Lemma 4.5 in \cite{KH01}), we also have $|B_n^*/B_n|\le 1$ in $\D$ which follows from the maximum principle and identity $|B_n^*/B_n|=1$ that holds on $\T$. This proves the lemma. \qed

\medskip

\begin{Lem}\label{l20}
If $X_n\subset \overline\D$ and
$
\lim_{n\to\infty}\sup_{z\in X_n}|F(z)- \widehat F_n(z)|= 0\,,
$
then
\[
\lim_{n\to\infty}\sup_{z\in X_n}|z(f(z) - \widehat f_{n}(z))|= 0.
\]
Conversely, if $\sup_{z\in \cup_{n\ge 1}X_n}|f(z)|<1$ and $\lim_{n\to\infty}\sup_{z\in X_n}|z(f(z)-\widehat f_n(z))|=0$, then
\[
\lim_{n\to\infty}\sup_{z\in X_n}|F(z)-\widehat F_n(z)|=0\,.
\]
\end{Lem}
\beginpf Recall that
\[
F=\frac{1+zf}{1-zf}, \qquad \widehat F_n=\frac{1+z\widehat f_n}{1-z\widehat f_n}.
\]
Thus, we have 
\[
|F-\widehat F_n|=\left|\frac{2z(f-\widehat f_n)}{(1-zf)(1-z\widehat f_n)}\right|\le \frac{2|z(f-\widehat f_n)|}{(1-|f|)(1-|\widehat f_n|)}.
\]
Analogously,
\[
|zf-z\widehat f_n|=\left|\frac{2(F-\widehat F_n)}{(1+F)(1+\widehat F_n)}\right|\le 2|F-\widehat F_n|,
\]
where we used the fact that $\Re F \ge 0$, $\Re \widehat F_n \ge 0$ in $\ov{\D}$. Now both claims are evident. \qed

\medskip

Later, we will need the following technical result.
\begin{Lem}\label{l103} 
Suppose function $G_n$ is analytic on $\mathfrak{D}_n=\{\eta: |\eta-in|<n\}$, continuous on $\overline{\mathfrak{D}}_n$, and $\Re G_n>0$ for every $n \ge 1$. Assume that there are constants $C_1$ and $C_2$ such that $\Re C_1>0$,
\begin{equation}\label{sd_oo}
\lim_{n\to\infty} G_n(\eta)=C_1
\end{equation}
 uniformly over compacts in $\C^+$,
\begin{equation}\label{sd_ww}
\lim_{n\to\infty}\arg G_n(in-ine^{it/n})= C_2, \quad 
\lim_{n\to\infty}\Bigl(\arg G_n(in-ine^{it/n})\Bigr)'= 0,
\end{equation}
and these two limits are uniform in $t$ over compacts in $\R$. Then, $C_2=\arg C_1$ and 
\begin{equation}\label{sd_kl}
\lim_{n\to\infty}\sup_{\eta\in \overline{\mathfrak{H}}_{b,n}}|G_n(\eta)-C_1|=0,
\end{equation}
for every $b>0$, where $\mathfrak{H}_{b,n}=\mathfrak{D}_n\cap \{\eta:|\eta|<b\}$.  
\end{Lem}
\beginpf The function $u_n=\Im \log G_n=\arg G_n$ is harmonic in $\mathfrak{D}_n$, continuous on $\overline{\mathfrak{D}}_n$, and $|u_n|\le \pi/2$. For every point $\eta\in \mathfrak{D}_n$, we can write Poisson formula
\[
u_n(\eta)=\int_{\partial \mathfrak{D}_n}u_n(\xi)d\omega_\eta(\xi),
\]
where $\omega_\eta$ is harmonic measure at $\eta$ for $\mathfrak{D}_n$ (the rescaled unit disk). The first condition in \eqref{sd_ww} and $|u_n|\le \pi/2$ imply that
\begin{equation}\label{sd_vv}
\lim_{n\to\infty}\sup_{\eta\in \overline{\mathfrak{H}}_{b,n}}|u_n(\eta)-C_2|=0,
\end{equation}
for every $b$. Thus, $C_2=\arg C_1$. Next, we will use the fact that the function analytic on the compact simply connected domain in $\C$ can be recovered from the boundary value of its imaginary part (up to a real constant). Indeed, let $\lambda_n(k)$ be conformal map of $\D$ to $\mathfrak{H}_{b,n}$ such that the lower arc of $\partial\mathfrak{H}_{b,n}$, i.e., points $\eta:\eta\in \partial\mathfrak{H}_{b,n}$ that satisfy $|\eta-in|=n$, corresponds to lower semicircle of $\partial \D$, i.e., $k: k\in \partial \D$ for which $\Im k<0$. Consider $\Gamma_n(k)=\log G_n(\lambda_n(k))$. It is analytic in $\D$, continuous in $\overline{\D}$ and, given conditions of the lemma, satisfies
\begin{equation}\label{sad_2019_1}
\lim_{n\to\infty}\Im\Gamma_n(e^{i\theta})=\arg C_1, \quad \lim_{n\to\infty}(\Im\Gamma_n(e^{i\theta}))'= 0\,,
\end{equation}
uniformly in $\theta\in [-\pi+\delta,-\delta]\cup[\delta,\pi-\delta]$ for every $\delta>0$. Moreover, $\lim_{n\to\infty}\Gamma_n(k)=\log C_1$ uniformly on compacts in $\D$ and $|\Im \Gamma_n|\le \pi/2$ in $\overline{\D}$. We can recover $\Gamma_n$ by the boundary values of its imaginary part as follows:
\[
\Gamma_n(k)=i\int_{\T} \Im \Gamma_n(\xi)\cdot\frac{1+\bar\xi k}{1-\bar\xi k}\,dm(\xi)+\Re \Gamma_n(0)\,.
\]
The conditions on $\{\Gamma_n\}$ and simple estimates on the integral above imply that $\lim_{n\to\infty}\Gamma_n(k)=\log C_1$ uniformly in $k: \{k\in \overline{\D}, |k-1|\ge \delta, |k+1|\ge \delta\}$, where $\delta$ is any positive number. In particular, this yields 
\[
\lim_{n\to\infty}\sup_{\eta\in \overline{\mathfrak{H}}_{b-1,n}}|G_n(\eta)-C_1|=0
\]
in the variable $\eta$.
Since $b$ is arbitrary positive, the lemma is proved.
\medskip

\begin{Lem}\label{sd_pp}
In Theorem \ref{t2}, if $(a)$, $(b)$, or $(c)$ holds, then $\lim_{n\to\infty}f_n(\xi) = 0$ for almost every $\xi\in \T$.
\end{Lem}
\beginpf From Lemmas \ref{lp1}, \ref{sd_33}, and \ref{sd_44}, we know that conditions $(a)$--$(c)$ are equivalent to each other. 
As before, let $\gamma_n$ denote the argument of the Blaschke product $b_n = \phi_n/\phi_n^*$ and let $\widetilde \gamma_n$ be the argument of $\widetilde b_n = \psi_n/\psi_n^*$. The proof of Lemma~\ref{sd_33} gives control for the derivatives of $\gamma_n$ and $\widetilde\gamma_n$ at almost every point $\xi\in \T$. Without loss of generality, assume that this point $\xi$ is equal to $1$. Additionally,  assume that the nontangential limit $f(1) = \lim_{z\to 1}f(z)$ exists and $|f(1)|<1$. This last condition implies existence of nontangential limit of $F$ at point $1$ and an estimate $\Re F(1)>0$. From the proof of Lemma \ref{sd_33}, we have 
\begin{equation}\label{sd_he1}
\lim_{n\to\infty}\arg \gamma_n'(\tau/n)/n = 1,\qquad \lim_{n\to\infty}\widetilde \gamma_n'(\tau/n)/n = 1
\end{equation}
 uniformly over compacts in $\R$. 

\medskip

By Lemma~\ref{lp1} and Lemma~\ref{l11}, we get $\lim_{n\to\infty}\sup_{z\in \Upsilon_{n^{-1},\rho,\alpha,\beta}}|f_n(z)|=0$ for arbitrary $\rho\in (0,1)$ and $0<\alpha<\beta$. From Lemma \ref{sd_55} and Lemma \ref{l20}, one has $\lim_{n\to\infty}\sup_{z\in \Upsilon_{n^{-1},\rho,\alpha,\beta}}|F(z)-\widehat F_n(z)|=0$. For $n \ge 1$ and $\eta$ in the disk $\mathfrak{D}_n$ defined in Lemma \ref{l103}, we set $G_n(\eta)=\widehat F_n(1-\eta/(in))$. The existence of nontangential limit of $F$ at point $z=1$ gives $\lim_{n\to\infty}G_n(\eta)=F(1)$ for every $\eta\in \C^+$. This convergence is in fact uniform over compacts in $\C^+$ by Lemma \ref{l11}. We will apply Lemma \ref{l103} next.
Notice that $\Re G_n>0$ in $\mathfrak{D}_n$. If we define 
$$
u_n(t)= \arg G_n(in - ine^{it/n})= \arg \widehat F_n(e^{it/n})=\arg \psi_n^*(e^{it/n})-\arg \phi_n^*(e^{it/n}),
$$ 
then $\lim_{n\to\infty}u_n'= 0$ uniformly over any compact in $\R$ as follows from \eqref{sd_he1} and a simple relation between the arguments of $b_n=\phi_n/\phi_n^*$ and $\phi_n^*$. Since $|u_n|\le \pi/2$, we can choose a subsequence $\{u_{n_j}\}$ such that 
\begin{equation}\label{sd_77}
\lim_{j\to\infty}u_{n_j} = C_*, \qquad \lim_{j\to\infty} u_{n_j}'= 0
\end{equation}
where $C$ is some constant and the both convergences are uniform on compact subsets of $\R$. Now apply Lemma \ref{l103} to $G_{n_j}$ taking $C_1=F(1)$ and $C_2=C_*$ to get 
\begin{equation}
\lim_{j\to\infty}\sup_{\eta\in \overline{\mathfrak{H}}_{b,n}}|{G}_{n_j}(\eta)-F(1)|=0
\end{equation}
for every $b$, where the sets $\mathfrak{H}_{b,n}$ are defined in Lemma \ref{l103}. Arguing by contradiction, we can strengthen this to 
\begin{equation}
\lim_{n\to\infty}\sup_{\eta\in \overline{\mathfrak{H}}_{b,n}}|{G}_{n}(\eta)-F(1)|=0
\end{equation}
for every $b$. Taking $\eta=0$, we get $\lim_{n\to\infty}\widehat F_n(1)= F(1)$ and thus $\lim_{n\to\infty} \widehat f_n(1) = f(1)$ by Lemma \ref{l20}. The last property is equivalent to $\lim_{n\to\infty} f_n(1) = 0$ by Lemma~4.8 in \cite{KH01}. \qed

\medskip
In the next lemma, we will control nontangential maximal function.\smallskip
\begin{Lem}\label{sd_pip} Let $\rho\in (0,1)$. If $g_n$ is analytic in $\D$, $|g_n|\le 1$, and $\lim\limits_{n\to\infty}g_n(\xi)=0$ for almost every $\xi\in \T$, then $\sup_{z\in S_\rho^*(\xi)}|g_n(z)|\to 0$ for almost every $\xi\in \T$. 
\end{Lem}
\beginpf By Egorov's theorem, for every $j\in \mathbb{N}$, we can find $E_j\subset \T$, $|E_j^c|<1/j$ and $\lim_{n\to\infty}g_n = 0$ on $E_j$ uniformly. We can assume without loss of generality that each point of $E_j$ is a Lebesgue point. Take $\xi\in E_j$ and let $z\in S_\rho^*(\xi)$. Write Poisson formula for harmonic function $g_n$:
$
g_n(z)=\P(g_n, z).
$
By dominated convergence theorem, we have $\lim_{n\to\infty} g_n = 0$ uniformly on compacts in $\D$. Then, 
\[
\sup_{|z|>r, z\in S_\rho^*(\xi)} |g_n(z)|\lesssim \delta(r)+\sup_{\eta\in E_j}|g_n(\eta)|\,,
\]
where $\lim_{r\to 1}\delta(r)= 0$ because $\xi$ is a Lebesgue point for $E_j$. For the second term, we have $\lim_{n\to\infty}\sup_{\eta\in E_j}|g_n(\eta)|=0$ due to uniform convergence to zero on $E_j$. Thus, 
\[
\sup_{z\in S_\rho^*(\xi)} |g_n(z)|=\sup_{|z|<r, z\in S_\rho^*(\xi)} |g_n(z)|+\sup_{|z|>r, z\in S_\rho^*(\xi)} |g_n(z)|
\]
and $\lim_{n\to\infty}\sup_{z\in S_\rho^*(\xi)} |g_n(z)|=0$
if we first fix $r$ close enough to $1$ and then let $n\to\infty$. Since $j$ is arbitrary, we get statement of the lemma.\qed

\medskip
Finally, we are ready to prove Theorem \ref{t2}.\smallskip

\noindent{\bf Proof of Theorem \ref{t2}}. Lemma \ref{lp1} shows that $(a)$ implies $(c)$. Lemma \ref{sd_33} shows that $(c)$ implies $(a)$ and $(b)$. Lemma \ref{sd_44} proves that $(b)$ implies $(a)$. This establishes equivalence of $(a)$, $(b)$, and $(c)$. Lemma~\ref{sd_117} shows that $(d)$ yields $(a)$. Finally, Lemmas~\ref{sd_pp} and \ref{sd_pip} prove that $(a)$, $(b)$, $(c)$ give $(d)$.
\qed



\bibliographystyle{plain} 
\bibliography{bibfile}
\enddocument